%
%
%
%
%
%
%
%
%
\scrollmode
\magnification=\magstep1
\parskip=\smallskipamount

\def\demo#1:{\par\medskip\noindent\it{#1}. \rm}
\def\ni{\noindent}               
\def\ll{\leftline}

\def\begin{\ll{}\vskip 10mm \nopagenumbers}  
\def\pn{\footline={\hss\tenrm\folio\hss}}   
\def\ii#1{\itemitem{#1}}

%
%
\outer\def\beginsection#1\par{\bigskip
  \message{#1}\leftline{\bf\&#1}
  \nobreak\smallskip\vskip-\parskip\noindent}

%
%
\outer\def\proclaim#1:#2\par{\medbreak\vskip-\parskip
    \noindent{\bf#1.\enspace}{\sl#2}
  \ifdim\lastskip<\medskipamount \removelastskip\penalty55\medskip\fi}
\def\endpr{\hfill $\spadesuit$ \medskip}

%
%

\def\R{{\rm I\kern-0.2em R\kern0.2em \kern-0.2em}}
\def\N{{\rm I\kern-0.2em N\kern0.2em \kern-0.2em}}
\def\P{{\rm I\kern-0.2em P\kern0.2em \kern-0.2em}}
\def\B{{\rm I\kern-0.2em B\kern0.2em \kern-0.2em}}
\def\C{{\rm C\kern-.4em {\vrule height1.4ex width.08em depth-.04ex}\;}}
\def\CP{\C\P}
\def\Z{{\bf Z}}

%
%
%
%
\def\cA{{\cal A}}

\def\cC{{\cal C}}

\def\cH{{\cal H}}

\def\cJ{{\cal J}}
\def\cK{{\cal K}}

\def\cO{{\cal O}}

\def\cU{{\cal U}}
\def\cV{{\cal V}}

%
%
%

\def\e{\epsilon}
\def\z{\zeta}

\def\l{\lambda}

%
%
%
%
\def\bar{\overline}              
\def\bs{\backslash}              
\def\di{\partial}                
\def\dibar{\bar\partial}         
\def\hra{\hookrightarrow}

%
%
\def\disc{\triangle}             

%
%
\def\dim{{\rm dim}\,}                    
\def\holo{holomorphic}                   
\def\nbd{neighborhood}                   
\def\psc{pseudoconvex}                   
\def\ra{real-analytic}                   
\def\spsh{strongly\ plurisubharmonic}
\def\hc#1{${\cH}({#1})$-convex}          
\def\ss{\subset\!\subset}                

\def\supp{{\rm supp}\,}                  

\def\iff{if and only if}
\def\hvf{holomorphic vector field}
\def\hvb{holomorphic vector bundle}

\def\phe{proper holomorphic embedding}

\def\wh{\widehat}
\def\wt{\widetilde}

\begin
\ll{\bf Oka's principle for holomorphic submersions}
\ll{\bf with sprays}

\bigskip\medskip
\ll{\bf Franc Forstneri\v c $\bf \cdotp$ Jasna Prezelj}
\medskip
\ll{\it Mathematics Subject Classification (2000): \rm 32H05, 32L05, 32Q28, 32Q55}
\medskip\bigskip

%
%
%
%
\ll{\bf Introduction}
\medskip

\ni The {\it Oka principle} is a fundamental principle
of complex analysis which says that on Stein manifolds 
(closed complex submanifolds of affine spaces), analytic problems 
of cohomological (or even homotopical) nature have only topological 
obstructions. The story began in 1939 when K.\ Oka proved that on a 
domain of holomorphy a second Cousin problem is solvable by holomorphic 
functions if it is solvable by continuous functions [Oka]. Oka's result
has the following equivalent formulation: 
{\it If $E\to X$ is a \holo\ $\C^*$-bundle over a domain of holomorphy 
then every continuous section of $E$ is homotopic to a \holo\ section.}
In 1957 H.\ Grauert [Gra] extended Oka's result to sections 
of \holo\ fiber bundles with complex homogeneous fibers over 
Stein manifolds (see also the papers [Car, FR, HL]).

In 1989 M.\ Gromov [Gr1] announced a major generalization
of Grauert's theorem, replacing complex Lie groups (and homogeneous
spaces) by a larger class of complex manifolds with a spray. 
A {\it spray} on a complex manifold $Y$ is a holomorphic vector bundle 
$p\colon E\to Y$ with a holomorphic map $s\colon E\to Y$ (the spray map)
which is the identity on the zero section and whose restriction to each fibre 
of $E$ is a submersion at zero. Such sprays exist on complex homogeneous 
spaces, on complements of affine algebraic subvarieties of codimension at 
least two and, more generally, on manifolds whose holomorphic tangent bundle 
is spanned by finitely many $\C$-complete holomorphic vector fields. 
The notion of spray extends in a natural way to submersions. 
The following is Theorem 4.5 in [Gr1]:
{\it If $h\colon Z\to X$ is a \holo\ submersion onto a Stein 
manifold $X$ such that every point $x\in X$ has an open \nbd\ 
$U\subset X$ such that $Z|_U=h^{-1}(U) \to U$ admits 
a fiber-dominating spray then the inclusion of the space of
holomorphic sections of $h$ to the space of continuous sections 
of $h$ is a weak homotopy equivalence.}  This will be referred 
to as the {\it parametric Oka principle} for sections of $Z\to X$.

In this paper we give a complete proof of this result and some 
extensions (Theorems 1.2 and 1.5). We also give examples 
(Theorem 1.6) and prove interpolation results 
(Theorems 1.7 and 1.9) which have been used in the construction 
of proper holomorphic embeddings of Stein manifolds into Euclidean 
spaces of minimal dimension [EGr, Sch, Pre1]. This paper is a sequel 
to [FP1] where we proved the same results for sections of fiber bundles 
whose fiber admits a dominating spray (thereby clarifying the 
results announced in sect.\ 2.9 of [Gr1]). The relevant analytic 
tools were proved in [FP1] following [Gra] and [Gr1]; they include a 
homotopy version of the Oka-Weil theorem for sections of submersions 
with sprays [FP1, sect.\ 4] and a gluing lemma for holomorphic 
sections over Cartan pairs [FP1, sect.\ 5]. 

To prove the Oka principle we develop an inductive scheme for patching 
collections of local holomorphic sections into semi-global \holo\ sections. 
Such a scheme was proposed in [Gr1, sect.\ 4], but we felt that an 
explanation would be welcome. This procedure requires that there
exist holomorphic homotopies between the individual sections in the given 
family, two-parameter homotopies between the one-parameter 
homotopies, etc., in order to insure that all `triangles of homotopies' 
which appear at various steps of the construction are `contractible'. 

To keep track of the data we introduce in section 3 the notion 
of holomorphic (resp.\ continuous) complexes and prisms. 
A {\it holomorphic complex} is a family of local holomorphic
sections of $h\colon Z\to X$, parametrized by the
{\it nerve} of an open covering $\cU=\{U_j\}$ of $X$.
A point in the parameter space determines an open set
in $X$ (the intersection of a certain collection of sets
from the covering $\cU$) on which the correponding section 
of $h$ is defined, and we have natural restriction conditions.
A global section $f$ is the same as a `constant complex'
in which any section is the restriction of $f$ to the appropriate 
set. A {\it holomorphic $k$-prism} is a homotopy of holomorphic complexes
with parameter in the $k$-dimensional cube. Similarly one defines
continuous complexes.

\pn

The procedure runs as follows. We begin by deforming the initial 
continuous section $a=a_{*,0}$ by a homotopy of continuous 
complexes $a_{*,t}$ ($t\in [0,1]$) to a holomorphic complex 
$a_{*,1}$ (Proposition 4.7). We then inductively 
construct a sequence of \holo\ complexes $f_*^n$ ($n\in \N$),
with $f_*^1=a_{*,1}$, such that $f_*^n$ is constant over the union 
of the first $n$ sets in the given covering of $X$
(i.e., it represents a \holo\ section there). When $n\to\infty$ 
the sequence $f_*^n$ converges uniformly on compacts in $X$ to a 
global \holo\ section $f\colon X\to Z$ of $h\colon Z\to X$.  
The heart of the proof is Proposition 5.1. 
This modification process requires special coverings of $X$, 
called {\it Cartan strings}, which were constructed 
by Henkin and Leiterer in [HL].

Grauert's Oka principle had numerous applications. Gromov's 
extensions have already been used in the embedding theorem for 
Stein manifolds into Euclidean spaces of minimal dimension [EGr, Sch]. 
Further references to recent applications are included at the end 
of the paper. Complete results on the Oka principle for 
maps of Riemann surfaces were obtained by Winkelmann [Win].

The Oka principle is a special case of the {\it homotopy principle}
whose validity for an analytic problem means that {\it an analytic 
solution exists provided that there are no topological obstructions}. 
Classical examples include the Smale--Hirsch theory of immersions 
of real manifolds and the theory of totally real and lagrangian 
immersions due to Lees and Gromov. A good reference are  the 
monographs [Gr2] and [Spr].

%
%
%
%

\beginsection 1. The results.

Let $h\colon Z \to X$ be a \holo\ submersion onto a Stein
manifold $X$. For each $x\in X$ we denote by $Z_x=h^{-1}(x)$
the fiber over $x$. At each point $z\in Z$ the tangent space
$T_z Z$ contains the {\it vertical tangent space}
$VT_z(Z) = \{ e\in T_z Z\colon Dh(z)e =0\} = T_z Z_{h(z)}$.         
If $p\colon E\to Z$ is a \hvb\ over $Z$, we denote by
$E_z=p^{-1}(z) \subset E$ its fiber over $z\in Z$ and
by $0_z \in E_z$ the zero element of $E_z$.

%
%
%
%
\medskip\ni \bf 1.1 Definition: \rm (Gromov [Gr1])  \sl
Let  $h\colon Z\to X$ be a \holo\ submersion of a complex manifold $Z$
onto a complex manifold $X$. A {\rm spray} on $Z$ associated to $h$
(or a fiber-spray) is a triple $(E,p,s)$, where $p\colon E\to Z$
is a \hvb\ and $s\colon E\to Z$ is a \holo\ map such that for
each $z\in Z$ we have
\item{(i)}   $s(E_z) \subset Z_{h(z)}$
(equivalently, $h\circ p=h\circ s$),
\item{(ii)}  $s(0_z)=z$, and
\item{(iii)} the derivative $Ds(0_z) \colon T_{0_z} E \to T_z Z$
maps the subspace $E_z \subset T_{0_z}E$ surjectively
onto the vertical tangent space $VT_z(Z)$.
\medskip\rm

\ni The restriction $VD s(z) = Ds(0_z)|E_z \colon E_z\to VT_z(Z)$
is called the {\it vertical derivative} of $s$ at $z\in Z$. 
Gromov's definition of a spray only includes properties (i) and (ii), 
and a spray which also satisfies the domination property (iii) 
is called a (fiber-) {\it dominating spray}.

%
%
%
%
\proclaim 1.2 Theorem:
\rm (Gromov [Gr1], 4.5 Main Theorem) \sl
Let $h\colon Z\to X$ be a \holo\ submersion of a complex manifold
$Z$ onto a Stein manifold $X$. Assume that each point in $X$ has an 
open \nbd\ $U\subset X$ such that the submersion $h^{-1}(U)\to U$
admits a (fiber dominating) spray. Then the inclusion
$\iota \colon {\rm Holo}(X,Z) \hra {\rm Cont}(X,Z)$
of the spaces of holomorphic sections into the space of continuous 
sections is a weak homotopy equivalence.

This means that $\iota$ induces an isomorphism of the respective
homotopy groups of the two spaces (which are endowed with the 
compact--open topology). In particular their path connected components 
are in one-to-one correspondence, which means that
\item{(i)}  any continuous section of $Z\to X$ can be
homotopically deformed to a \holo\ section, and
\item{(ii)} any homotopy of sections $f_t\colon X\to Z$ ($0\le t\le 1$)
between two \holo\ sections $f_0$ and $f_1$ can be deformed into
another homotopy consisting of holomorphic sections.
\rm\medskip

In the special case when $Z\to X$ is a fiber bundle whose fiber 
admits a spray, a complete proof of Theorem 1.2 was given in [FP1] 
(thereby clarifying the results announced in sect.\ 2.9 of [Gr1]). 
When the conclusion of Theorem 1.2 holds, we shall say 
that sections of $h$ satisfy the {\it parametric Oka principle}. 
Likewise we say that the (parametric) Oka principle holds for maps $X\to Y$ 
of a Stein manifold $X$ into a complex manifold $Y$ if it holds 
for sections of the trivial fibration $Z=X\times Y\to X$.

Examples of spaces with sprays can be found in [Gr1, section 4.6.B]
and in [FP1]. We recall the following example.

%
%
\demo Example 1: Let $h\colon Z\to X$ be a \holo\ submersion.
Suppose that $Z$ admits finitely many $\C$-complete \hvf s
$V_1,V_2,\cdots,V_N$ which are vertical (tangent to $VT(Z)$) and
which span $VT_z(Z)$ at each point $z\in Z$. $\C$-completeness
means that the flow $\phi^t_j$ of $V_j$ is defined for all
complex values of the time parameter $t$.
The map $s\colon Z\times \C^N \to Z$, defined by
$$
    s(z;t_1,\ldots,t_N) =
      \phi_1^{t_1}\circ \phi_2^{t_2}
        \circ\cdots\circ \phi_N^{t_N}(z),
$$
satisfies $s(z;0,\ldots,0)=z$ and
${\di\over \di t_j} s(z;0,\ldots,0)=V_j(z)$
for $z\in Z$ and $1\le j\le N$. Since these vectors
span $VT_z(Z)$, $s$ is a spray on $Z$.
\endpr

%
%
\proclaim 1.3 Corollary:
Let $h\colon Z\to X$ be a \holo\ submersion onto a Stein
manifold $X$. Assume that each point in $X$ has an open
\nbd\ $U\subset X$ such that there exist finitely many $\C$-complete
\hvf s on $Z|U$ which are vertical with respect to $h$ and
which span the vertical tangent space $VT_z(Z)$ at each
point $z\in Z|U$. Then the conclusion of Theorem 1.2 holds.

We now describe a more general version of the Oka principle
(compare with Theorem 1.4 and Corollary 1.5 in [FP1]).
Let $h\colon Z\to X$ be a holomorphic submersion and let $P$
be a compact Hausdorff space (the parameter space).
Our basic objects now are continuous maps $f\colon X\times  P\to Z$
such that $f(\cdotp,p)\colon X\to Z$ is a section of $h$ for
each fixed $p\in P$. A {\it homotopy} of such maps is a
continuous map $H\colon X\times P\times [0,1]\to Z$
such that $H_t(\cdotp,p)=H(\cdotp,p,t) \colon X\to Z$ is
a section of $h$ for all $p\in P$ and $t\in [0,1]$.

We denote by ${\cH}(X)$ the algebra of all \holo\ functions on
$X$. A compact set $K \subset X$ is called  \hc{X} 
(or {\it holomorphically convex}) if for each $x\in X\bs K$ there 
is an $f \in {\cH}(X)$ with $|f(x)|>\sup_K |f|$.

%
%

\proclaim 1.4 Definition: A subset $P_0$ in a topological space
$P$ is called {\rm nice} if there exists an open set $U\subset P$
containing $P_0$ and a strong deformation retraction of $U$
onto $P_0$. The empty subset in $P$ will be considered nice.

%
%
%
%
\proclaim 1.5 Theorem:  Let $X$ be a Stein manifold, $K\subset X$
a compact \hc{X}\ subset and $h\colon Z\to X$ a \holo\ submersion
onto $X$. Assume that for each point $x\in X\bs K$ there is a
\nbd\ $U_x\subset X$ such that the submersion
$h\colon h^{-1}(U_x) \to U_x$ admits a spray (def.\ 1.1).
Let $P$ be a compact Hausdorff space and $P_0 \subset P$ a nice
compact subset (def.\ 1.4). Assume that $a\colon X\times P\to Z$
is a continuous map such that for each $p\in P$,
$a(\cdotp,p) \colon X\to Z$ is a
section of $h\colon Z\to X$ which is \holo\ in an open
set $U_0\supset K$, and
$a(\cdotp,p)$ is \holo\ on $X$ for each $p\in P_0$.
Let $d$ be a metric on $Z$ compatible with the manifold
topology. Then for each $\e>0$ there exists a homotopy
$H_t\colon X\times P\to Z$ $(t\in [0,1])$ satisfying:
\item{(i)}   $H_0=a$,
\item{(ii)}  the section $H_1(\cdotp,p)\colon X\to Z$ is
\holo\ for each $p\in P$,
\item{(iii)} the homotopy is fixed on $P_0$, i.e., $H_t(\cdotp,p)$
is independent of $t$ for $p\in P_0$, and
\item{(iv)}  $d\bigl( H_t(x,p), a(x,p) \bigr) < \e$ for all
$x\in K$, $p\in P$ and $0\le t\le 1$.

\demo Theorem 1.5 implies Theorem 1.2:
If we take $P$ to be the $n$-sphere $S^n$ and $P_0=\emptyset$,
Theorem 1.5 shows that each continuous map $S^n\to {\rm Cont}(X,Z)$
can be homotopically deformed to a map $S^n\to {\rm Holo}(X,Z)$.
Similarly, if $P$ is the closed $(n+1)$-ball $B^{n+1}\subset\R^{n+1}$
and $P_0=\di B^{n+1}=S^n$, Theorem 1.5 shows that each map
$S^n\to {\rm Holo}(X,Z)$ which extends to a map
$B^{n+1}\to {\rm Cont}(X,Z)$ also extends to a map
$B^{n+1}\to {\rm Holo}(X,Z)$. This means that the inclusion
${\rm Holo}(X,Z) \hra {\rm Cont}(X,Z)$ induces an isomorphism 
of the respective homotopy groups of the two spaces as 
claimed by Theorem 1.2.
\endpr

Next we consider the validity of the Oka principle for maps
of Stein manifolds into $\C^q\bs\Sigma$, where $\Sigma$ is
a closed complex subvariety of $\C^q$. Immediate examples show that
the Oka principle fails in general when $\Sigma$ is a complex hypersurface,
except in special cases such as when a complex Lie group acts
transitively on $\C^q\bs\Sigma$. The situation is fairly complicated
even for subsets of higher codimension as the following result shows:

%
%
%
%
\medskip\ni\bf 1.6 Theorem. \sl
(a) If $q\ge2$ and $\Sigma\subset \C^q$ is a closed analytic
subset of complex codimension $\ge 2$ in $\C^q$ such that,
with respect to some holomorphic coordinates
$z=(z',z_q)\in \C^q$ and some constant $C>0$ we have
$$
    \Sigma \subset  \Gamma = \{z\in \C^q \colon
    |z_q|\le C(1+|z'|)\},                          \eqno(1.1)
$$
then the Oka principle holds for maps from any Stein manifold
into $\C^q\bs\Sigma$. This is the case in particular if
$\Sigma$ is an algebraic subset of codimension at least
two in $\C^q$.

\item{(b)} For each $q\ge 1$ there exist discrete sets
$\Sigma\subset \C^q$ such that the Oka principle fails
for maps of $X=(\C^*)^{2q-1}$ into $\C^q\bs \Sigma$.

\item{(c)} For each $1\le k< q$ there exist \phe s
$g\colon \C^k\hra \C^q$ such that the Oka principle fails
for maps of $X=(\C^*)^{2(q-k)-1}$ into $\C^q\bs g(\C^k)$.
\medskip\rm

\demo Remark: In the forthcoming paper [Fo5] we extend 
the Oka principle to maps from Stein manifolds into 
complements $\CP^n\bs \Sigma$ of projective-algebraic
subvarieties of codimension at least two.
\endpr

We now consider the problem of avoiding analytic subsets
in \hvb s by graphs of \holo\ sections. Theorems 1.7
and 1.9 below have been used in [EGr], [Sch] and [Pre1, Pre2].

%
%
%
%
\medskip\ni\bf 1.7 Theorem. \sl
Let $\pi\colon V\to X$ be a \hvb\ of rank $q\ge 2$ over a
Stein manifold $X$, let $X_0 \subset X$ be a closed complex subvariety
in $X$, and let $K \subset X$ be a compact \hc{X}\ subset.
Set $X'=X\bs (X_0\cup K)$ and $V'= \pi^{-1}(X') \subset V$.
Suppose that $\Sigma \subset V'$ is a closed analytic subset in $V'$
such that for each $x\in X'$, the fiber $\Sigma_x=\Sigma\cap V_x$
has complex codimension at least two in $V_x$, and there are an
open \nbd\ $U\subset X'$ of $x$, a set $\Gamma \subset \C^q$ as
in (1.1), and a biholomorphic map
$\Phi\colon U'= \pi^{-1}(U) \to U\times\C^q$
of the form $\Phi(x,v)=(x,\phi(x,v))$ ($x\in U,\ v\in V_x$),
satisfying
$$
    \Phi(U' \cap \Sigma) \subset U\times \Gamma.       \eqno(1.2)
$$
Given a continuous section $f_0\colon X\to V\bs \bar \Sigma$ that is \holo\
in an open set $U_0\supset X_0\cup K$, there is for each $k\in \Z_+$
a homotopy of continuous sections $f_t\colon X\to V\bs \bar\Sigma$
($0\le t\le 1$) which are \holo\ in a \nbd\ of $X_0\cup K$, which match
$f_0$ to order $k$ along $X_0$ and approximate $f_0$ uniformly on $K$,
and the section $f_1$ is \holo\ on $X$.
\medskip\rm

We emphasize that $\Sigma$ is not assumed to be closed in $V$,
and its closure $\bar \Sigma$ need not be an analytic subset of $V$.
We do not assume anything about the projection $\pi(\Sigma)$,
and the sets $V_x\bs \Sigma_x$ need not be biholomorphically
equivalent for different base points $x$. The condition (1.2)
will insure the existence of a spray on $U'\bs \Sigma$.

We state an important special case of Theorem 1.7. To each
\hvb\ $\pi\colon V\to X$ of rank $q$ we can associate a holomorphic
fiber bundle $\bar\pi\colon \bar V\to X$ whose fiber
$\bar V_x \cong \CP^q$ is obtained by compactifying $V_x\cong\C^q$
with the hyperplane at infinity. Since every complex linear automorphism
of $\C^q$ extends to a unique projective automorphism of $\CP^q$,
the fibers $\bar V_x$ patch together into a \holo\ bundle over
$X$.

\proclaim 1.8 Corollary: Let $\pi\colon V\to X$ be a \hvb\
of rank $q$ over a Stein manifold $X$, and let
$\bar\pi\colon \bar V\to X$ be the associated bundle with fiber $\CP^q$
as above. If $\Sigma\subset \bar V$ is a closed analytic subset such
that for each $x\in X$ the fiber $\Sigma_x=\Sigma\cap \bar V_x$ is
of complex codimension at least two in $\bar V_x$ then the Oka principle
holds for sections $f\colon X\to V\bs \Sigma$ (in the strong form
described by Theorem 1.7).

We will show that, if $\Sigma$ is as in the corollary, the condition (1.2)
holds with respect to some local vector bundle charts $\Phi$ on $V$ and
hence we get sprays. By Chow's theorem [Chi] each fiber $\Sigma_x$ is an
algebraic subset of $\bar V_x$, and the codimension condition insures
that $\Sigma_x$ does not contain the hyperplane at infinity.
This last condition can be equivalently described as follows.
In a local bundle chart on $\pi^{-1}(U)\cong U\times \C^q$ over a
small subset $U\subset X$, the set $\Sigma\cap \pi^{-1}(U)$ is given
by finitely many equations
$$
    g_j(x,w)= \sum_{|I|\le d_j} g_{j,I}(x) w^I =0
    \quad (1\le j\le m)
$$
which are polynomial in $w \in \C^q$ and where the coefficients
$g_{j,I}$ are holomorphic in $U$.
The closure $\bar \Sigma \subset U\times\CP^q$ is defined by
the corresponding homogenized equations, obtained by adding a suitable
power of an additional variable $w_0$ to each term (so $w_0=0$
is the hyperplane at infinity). Let $P_j$ be the top order homogeneous
part of $g_j$ (with respect to $w$); this is the part of $g_j$
which gets not $w_0$ terms after homogenization.
We see that $\bar \Sigma_x \subset \CP^q$
contains the hyperplane at infinity $\{w_0=0\}$ \iff\ $P_j(x,w)=0$ for
all $w\in \C^q$ and $j=1,\ldots, m$. The hypothesis precludes this 
from happening and hence Lemma 7.1 applies.

\demo Example 2: If $D=\{z\in \C\colon |z|<1\}$ and
$\Sigma\subset D\times \CP^1$ is a closed analytic subset
of dimension one (a complex curve), it is easy to show that
$\Sigma$ can be avoided by graphs of continuous (or smooth)
functions $f\colon D\to \C$, but in general there is no
\holo\ map $f\colon D\to \CP^1$ (i.e., meromorphic function on $D$)
whose graph avoids $\Sigma$, due to hyperbolicity.
In fact, let $\Sigma_n$ be the union of graphs of the following
functions on $D$: $z\to 0$, $z\to 1$, $z\to\infty$,
and $z\to kz$ for $k=1,2,\ldots,n$. Then for
sufficiently large $n$ there is no \holo\ (or meromorphic)
function on $D$ whose graphs avoids $\Sigma_n$.
Proof: if the graph of $f_n$ avoids $\Sigma_n$ for each
$n$, the sequence $\{f_n\}$ is a normal family, and hence
a subsequence $f_{n_k}$ converges to a \holo\ function
$f\colon D\to \C$. For sufficiently large $n$ the
complex line $w=nz$ intersects the graph of $f$
transversely at some point; the same line then intersects the
graph of $f_l$ for all sufficiently large $l$, a contradiction.
\medskip

In our last theorem the hypotheses are similar as in Theorem 1.7,
except that we replace the codimension condition on the fibers
$\Sigma_x$ by a homogeneity condition (compare with sect.\ 1 in [EGr]).
We state it without approximation on a \hc{X}\ set, although
the result actually holds in the same form as Theorem 1.7.

\medskip\ni \bf  1.9 Theorem. \sl
Let $\pi\colon V\to X$ be a \hvb\ over a Stein manifold $X$,
let $X_0\subset X_1\subset X$ be closed complex subvarieties of $X$,
and let $\Sigma\subset V$ be a closed complex subvariety such that
$\pi(\Sigma)=X_1\bs X_0$. Assume that for each point
$x_0\in X_1\bs X_0$ there are an open \nbd\ $U\subset X$
and a \holo\ action of a complex Lie group $G$ on
$U'=\pi^{-1}(U) \subset V$, satisfying:
\item{(i)}  $\pi(g\cdotp v)=\pi(v)$ for each $v\in U'$
and $g\in G$, and
\item{(ii)} for each $x\in U\cap X_1$,
$G$ preserves $V_x\bs \Sigma_x$ and acts transitively there.

\ni Then for any continuous section  $f_0\colon X\to V\bs\Sigma$ that
is \holo\ in a \nbd\ of $X_0$ and for any integer $k\in Z_+$
there is a homotopy $f_t\colon X\to V\bs\Sigma$ ($0\le t\le 1$)
satisfying the conclusion of Theorem 1.7.
\rm\medskip

The paper is organized as follows. In section 2 we explain the
idea of the proof of Theorems 1.2 and 1.5. In section 3 we
introduce the concept of a \holo\ (resp.\ continuous) complex and
prism, and in sect.\ 4 we introduce the notion of a Cartan string.
The heart of the proof is Proposition 5.1 and the subsequent results
of section 5. We conclude the proof of Theorem 1.5 by an inductive
construction in section 6. In section 7 we prove Theorem 1.6, and
in section 8 we prove Theorems 1.7 and 1.9.

%
%
%
%
\beginsection 2. Outline of the proof of the main theorem.

In this section we explain the scheme of proof of Theorems
1.2 and 1.5. We concentrate on the case without parameters
(i.e., when $P$ is a single point). Thus, given a compact
\hc{X}\ set $K\subset X$ and a continuous section
$a\colon X\to Z$ which is \holo\ in an open set $U_0\supset K$,
we shall construct a homotopy of sections $H^s \colon X\to Z$ ($0\le s\le 1$)
such that $H^0=a$, the section $f=H^1$ is \holo\ on $X$, and every
section $H^s$ is \holo\ near $K$ and it approximates $a$
uniformly on $K$.

We begin by constructing a collection of \holo\ sections
$a_{(j)}\colon U_j\to Z$ ($j=0,1,2,\ldots$) on a suitably chosen
open covering $\cU=\{U_0,U_1,U_2,\ldots\}$ of $X$, where $U_0 \supset K$
is the given set on which $a$ is already \holo, such that
for each $j\in \Z_+$ there is a homotopy of continuous sections
$a_{(j),s}\colon U_j\to Z$ ($0\le s\le 1$) satisfying
$a_{(j),0}=a|U_j$ and $a_{(j),1}=a_{(j)}$. For $j=0$ we
take $a_{(0),s}=a|U_0$ for all $s$. (See Proposition 4.7.)

We think of the collection $\{a_{(j)}\colon j\in \Z_+ \}$ as a
puzzle whose pieces should be rearranged into a \holo\
section $f\colon X\to Z$ which is homotopic to the initial section $a$.

Let us first consider what is involved in patching a pair of local
\holo\ sections into a single \holo\ section over the union of the two sets.
To simplify the notation we assume that $A,B\subset X$ are
compact subsets, $U\supset A$ and $V\supset B$ are open sets,
and $a\colon U\to Z$, $b\colon V\to Z$ are \holo\ sections of
$Z\to X$ over $U$ resp.\ $V$. We wish to replace the pair $(a,b)$
by a section $\wt a$ which is \holo\ in a \nbd\ of $A\cup B$ and 
which approximates $a$ on $A$. (The problem is nontrivial 
only when $C=A\cap B\ne\emptyset$.) We proceed in two steps. 
First we replace $b$ by another \holo\ section
$b_1\colon V\to Z$ (perhaps shrinking $V$ around $B$ if necessary)
which approximates $a$ very closely in some \nbd\ $W$ of $A\cap B$. 
In the second step we `glue' $a$ with $b_1$ into a new \holo\ 
section $\wt a$ over $A\cup B$. The first step can be achieved if 
the following hold:
\item{--} the set $C=A\cap B$ is Runge in $B$, i.e., we can
approximate \holo\ functions in small \nbd s of $C$ by
functions \holo\ in a \nbd\ of $B$,
\item{--} the submersion $Z\to X$ admits a spray over a \nbd\ of
$B$, and
\item{--} there is a \nbd\ $W\supset A\cap B$ and a homotopy
of \holo\ sections $b_t\colon W\to Z$ connecting
$b_0=b|W$ and $b_1=a|W$.

\ni
Granted these conditions, the h-Runge theorem (see Theorems 4.1
and 4.2 in [FP1]) shows that we can approximate the homotopy
$\{b_t\}$ uniformly in a \nbd\ of $C=A\cap B$ by a homotopy of sections
$\wt b_t$ ($0\le t\le 1$) which are \holo\ in a \nbd\ of
$B$. In particular the section $\wt b_1$ approximates $a$ as
well as desired in some fixed \nbd\ of $C$. Replacing $b$
by $\wt b_1$ we may thus assume that $b$ approximates
$a$ near $C$.

To glue $a$ and $b$ we linearize the problem and solve a
certain $\dibar$-equation in a \nbd\ of $A\cup B$.
This is accomplished by Theorems 5.1 and 5.5 in [FP1],
provided that $(A,B)$ is a {\it Cartan pair} in $X$
(def.\ 4.1 below). This means that each of the sets $A,B$ and 
$A\cup B$ has a basis of Stein neighborhoods, $C$ is Runge in $B$, 
and $(\bar{A\bs B})\cap (\bar{B\bs A})=\emptyset$.

In order to glue the local \holo\ sections
$a_{(j)}\colon U_j\to Z$ into a single \holo\ section
$f\colon X\to Z$ we perform the above steps inductively.
At each step the sets $U_j$ on which our sections are defined
may shrink, and we must control the shrinking so that in the
end we still have a covering of $X$. For this reason we initially
choose a special locally finite covering $\cA=\{A_0,A_1,A_2,\ldots\}$
of $X$ by {\it compact} sets satisfying
\item{(i)}  $K\subset A_0 \subset U_0$;
\item{(ii)} for each $n\ge 1$, the ordered collection of
sets $(A_0,A_1,\ldots,A_n)$ is a {\it Cartan string}
(def.\ 4.2 below). This property enables us to carry out
the gluing procedure by induction on $n$. In particular this
means that for each $n\ge 1$, $A^{n-1}=A_0\cup A_1\cup \cdots\cup A_{n-1}$
and $A_n$ form a Cartan pair;
\item{(iii)} for each $j\in \Z_+$ there is a \holo\ section
$a_{(j)}$ in an open \nbd\ $U_j\supset A_j$ which is homotopic
to $a|U_j$, and for $j\ge 1$ the restriction $Z|U_j$ admits a spray;
\item{(iv)}  for each pair $i\ne j$ such that
$A_i\cap A_j\ne\emptyset$, there is a \holo\ homotopy
between $a_{(i)}$ and $a_{(j)}$ in $U_{(i,j)}=U_i\cap U_j$; etc.

In general, for each multiindex $J=(j_0,j_1,\ldots,j_n)$
such that $A_J=A_{j_0}\cap\cdots\cap A_{j_n}\ne \emptyset$,
there is an $n$-dimensional \holo\ homotopy $a_J(t)$ in
$U_J= U_{j_0}\cap\cdots\cap U_{j_n}$, with the parameter
$t$ belonging to the standard $n$-simplex $\disc^n\subset \R^n$,
such that for $t$ belonging to a boundary face of
$\disc^n$ determined by a shorter multiindex $I\subset J$
we have $a_J(t)=a_I(t)|U_J$. The parameter space of our entire
colection of \holo\ sections and homotopies between them is the
geometric realization of the simplicial complex called the
{\it nerve} of the covering $\cA$ (see sect.\ 3).
The sets $U_j$ will shrink but the $A_j$'s will
stay the same during the entire construction.

Suppose inductively that we have already joined the
sections $a_{(0)},\ldots,a_{(n-1)}$
into a \holo\ section $f^{n-1}$ in a \nbd\ of
$A^{n-1}=A_0\cup A_1\cup\cdots\cup A_{n-1}$, using the
homotopies between them and the gluing procedure. We emphasize
that all modifications are done by \holo\ homotopies.
In the next step we must glue $f^{n-1}$ with the section $a_{(n)}$.
For this we need a \holo\ homotopy between the two sections in a \nbd\
of $A^{n-1}\cap A_n$. In the special case when $h\colon Z\to X$ is 
a fiber bundle whose fiber admits a spray such a homotopy can be 
constructed from a continuous homotopy between the two sections,
provided that our sets are chosen sufficiently carefully
(i.e., $A_n$ must be a \psc\ bump on $A^{n-1}$). This was explained
in [FP1]; see especially Proposition 6.1 there. The argument
in [FP1] does not seem to carry over to the present situation and
we need an alternative systematic way of insuring the existence
of such homotopies.

Our inductive construction is such that for each
$j=0,1,\ldots,n-1$ we have a \holo\ homotopy  between
$f^{n-1}$ and $a_{(n)}$ in a \nbd\ of $A_j\cap A_n$,
inherited from the initial homotopy between $a_{(j)}$ and $a_{(n)}$.
We now patch these $n$ partial homotopies into a single homotopy
defined in a \nbd\ of
$A^{n-1}\cap A_n =\cup_{j=0}^{n-1} (A_j\cap A_n)$.
This can be done in the same way as above by induction on $n$,
provided that the ordered collection
$(A_0\cap A_n,A_1\cap A_n,\ldots,A_{n-1}\cap A_n)$ is also
a Cartan string. Finally, since $(A^{n-1},A_n)$ is
a Cartan pair, we can glue $f^{n-1}$ and $a_{(n)}$ into a section
$f^n$ in a \nbd\ of $A^n$ which approximates $f^{n-1}$ on $A^{n-1}$.
This completes the induction step. (The details are given in sect.\ 5.)
The sequence of sections $f^n$ obtained in this way converges
uniformly on compacts in $X$ to a \holo\ section
$f \colon X\to Z$ which solves the problem (sect.\ 6).

%
%
%
%
\beginsection 3. Holomorphic complexes and prisms.

%
%
\proclaim 3.1 Definition: Let $\cA=\{A_0,A_1,A_2,\ldots\}$ be
any finite or countable family of nonempty subsets of $X$ which
is locally finite. The {\rm nerve} of $\cA$ is the (combinatorial)
simplicial complex $\cK(\cA)$ consisting of precisely those
multiindices $J=(j_0,j_1,\ldots,j_k) \in \Z_+^{k+1}$
$(k\in \Z_+)$ with increasing entries $0\le j_0<j_1<\cdots < j_k$
for which
$$ A_J=A_{j_0}\cap A_{j_1}\cap\cdots\cap A_{j_k} \ne \emptyset. $$
We denote the geometric realization of $\cK(\cA)$ by $K(\cA)$.

%
%
For simplicial complexes and their realization we refer the reader
to [HW]; here we only recall a few basic ideas.
The geometric realization of a simplicial complex is
a topological space which is a union of topological simplexes
of various dimensions such that any two of them either do not
interesect, they intersect along a simplex of lower dimension,
or one is contained in the other.

In the nerve complex, each multiindex $J=(j_0,j_1,\ldots,j_k)\in \cK(\cA)$
of length $k+1$ determines a closed $k$-dimensional face
$|J|\subset K(\cA)$, called the {\it body} (or {\it carrier})
of $J$, and $J$ is called the {\it vertex scheme} of $|J|$.
$|J|$ is homeomorphic to the standard $k$-simplex
$\disc^k \subset \R^k$ (the closed convex hull of the set
$\{0,e_1,e_2,\ldots,e_k\} \subset\R^k$ containing the origin and
the standard basis vectors), and any $k$-dimensional face
of $K(\cA)$ is of this form. Thus the topological space $K(\cA)$ is the
body of the combinatorial object $\cK(\cA)$. The vertices
of $K(\cA)$ correspond to the individual sets $A_j\in \cA$,
i.e., to singletons $\{j\}\in \cK(\cA)$. Given $I,J\in \cK(\cA)$
we have $|I|\cap |J|=|I\cap J|$. Thus for any two (bodies of) simplexes
in $K(\cA)$, either one is a subset of the other or else their
intersection is a simplex of lower dimension, possibly empty.

For each $n\in \Z_+$ we denote by
$$
    \cK^n(\cA) = \cK(A_0,A_1,\ldots,A_n) \subset \cK(\cA)
$$
the nerve of the finite subcollection $\cA_n=\{A_0,\ldots,A_n\}$
and by $K^n(\cA)$ its body. Occasionally we delete ${\cA}$ in the
notation. Clearly $\cK^n\subset \cK^{n+1}$ for
each $n$, and $\cK(\cA)=\bigcup_{n=0}^\infty \cK^n$.
More generally, for any multiindex $J=(j_0,j_1,\ldots,j_k) \in \Z_+^{k+1}$
(not necessarily belonging to $\cK(\cA)$) we denote by
$$
    \cK_J(\cA)=\cK(A_{j_0},A_{j_1},\ldots,A_{j_k})
$$
the nerve of the indicated collection of sets and by $K_J(\cA)$
its body. Note that $\cK_J(\cA)$ is a finite subcomplex of
$\cK(\cA)$ whose body is a $k$-dimensional simplex \iff\
$J\in \cK(\cA)$; otherwise it is a union of simplexes of
lower dimension. Occasionally we shall write simply $\cK$
instead of $\cK(\cA)$ when it is clear from the context
which collection $\cA$ is meant.

From now on we assume that $\cA=\{A_0,A_1,A_2,\ldots\}$ is a
(finite or countable) locally finite family of compact subsets of $X$.
Later on, $\cA$ will be a covering of $X$ with
additional properties, but this is not important for the
moment. An open \nbd\  of $\cA$ is a collection
$\cU=\{U_0,U_1,U_2,\ldots\}$ of open subsets of $X$
(with the same index set) such that $A_i\subset U_i$ for
each $i$. Such an \nbd\ is called {\it faithful}
if $\cK(\cU)=\cK(\cA)$, that is, the two families have the
same nerve complex. Clearly any locally finite family $\cA$
as has an open faithful \nbd\ $\cU$. As before we write
for each $J=(j_0,j_1,\ldots,j_k)$
$$
   U^J=U_{j_0}\cup U_{j_1}\cup\cdots \cup U_{j_k},\quad
   U_J=U_{j_0}\cap U_{j_1}\cap\cdots \cap U_{j_k}.
$$
If $h\colon Z\to X$ is a holomorphic submersion onto $X$
and $U\subset X$ is an open subset, we denote by
$\cO_h(U,Z)$ (resp.\ $\cC_h(U,Z)$) the set of
all \holo\ (resp.\ continuous) sections $f\colon U\to Z$
of $h \colon Z\to X$ over $U$.

%
%
\medskip \ni\bf 3.2 Definition. \sl
Let $h\colon Z\to X$ be a holomorphic map of a
complex manifold $Z$ onto a complex manifold $X$ and
let $\cA=\{A_0,A_1,A_2,\ldots\}$ be a locally finite,
at most countable family of compact sets in $X$.

\item{(i)} A {\rm \holo\ $\cK(\cA)$-complex} 
with values in $Z$ is a collection 
$$
  f_* = \{f_J \colon |J| \to \cO_h(U_J,Z),\ \ J \in \cK(\cA) \},
$$
where $\cU=\{U_0,U_1,U_2,\ldots\}$ is a faithful \nbd\ of $\cA$,
satisfying the following compatibility conditions:
$$
    I,J\in \cK(\cA),\ I\subset J\Longrightarrow
      f_J(t) = f_I(t)|_{U_J}, \ t\in |I|.   
$$
A {\rm continuous $\cK(\cA)$-complex} with values in $Z$
is a collection 
$$
   f_* = \{f_J \colon |J| \to \cC_h(U_J,Z),\ \ J \in \cK(\cA) \}
$$
satisfying the same compatibility conditions.

\item{(ii)} If $f_*$ is a (\holo\ or continuous)
$\cK(\cA)$-complex and $\cK'\subset \cK(\cA)$ is a subcomplex
of $\cK(\cA)$, we denote by $f_*|\cK'$ the restriction of
$f_*$ to $\cK'$.

\item{(iii)} A $\cK(\cA)$-complex $f_*$, defined on a faithful
\nbd\ $\cU$ of $\cA$, is {\rm constant} if there is a
section $g\colon \cup\{U\colon U\in \cU\} \to Z$ such that
$f_J(t)|U_J = g|U_J$ for each $J \in \cK(\cA)$ and each $t\in |J|$.

\demo Remark: We shall consider $\cK(\cA)$-complexes
in the sense of their {\it germs} on the sets in $\cA$.
Thus, two $\cK(\cA)$-complexes $f_*$ and $g_*$ are considered
equivalent if there is an open faithful \nbd\ $\cU=\{U_i\}$
of $\cA=\{A_i\}$ such that for each $J\in \cK(\cA)$
and $t\in |J|$, the sections $f_J(t)$ and $g_J(t)$ are
defined and equal in $U_J$. In practice we shall
not distinguish between equivalent complexes.
\endpr

A \holo\ complex is a bookkeeping tool that we shall need
to keep track of the sections and homotopies between them alluded
to in section 2. For each $J=(j_0,\ldots,j_k)\in \cK(\cA)$
we have a family of \holo\ sections
$$
    f_J(t) \colon U_J\to Z, \quad t\in |J|,
$$
depending continuously on the parameter $t\in |J|$, which we
may think of as a homotopy of holomorphic sections over the set
$U_J=U_{j_0}\cap U_{j_1}\cap\cdots \cap U_{j_k}$, with the
parameter $t$ belonging to the simplex $|J| \subset K(\cA)$.
For each face $I\subset J$ of $J$ and for
$t\in |I| \subset |J|$ the section $f_J(t)$ agrees with the
section $f_I(t)$, restricted from its original domain
$U_I \subset X$ to the subdomain $U_J$. It is worthwhile
to consider separately the one dimensional case.

%
%
\demo  Example 3: Let $J=(j_0,j_1) \in \cK(\cA)$ be a simplex of
length two. Its body $|J| \subset K(\cA)$ is a segment which we
can represent by $[0,1]\subset \R$. $J$ contains two zero
dimensional faces, namely the vertices $(j_0)$ and $(j_1)$
(corresponding to the sets $A_{j_0}$ resp.\ $A_{j_1}$),
and the bodies of these faces are identified with the
endpoints $\{0\}$ resp.\ $\{1\}$ of $[0,1]$.
The restriction of a $\cK(\cA)$-complex $f_*$ to the
subcomplex $\cK_J(\cA)\subset \cK(\cA)$, determined
by $J=(j_0,j_1)$, consists of a one parameter family
(a homotopy) $f_J(t)$ ($t\in |J|=[0,1]$) of sections over
$U_J= U_{j_0} \cap U_{j_1}$ such that $f_J(0)$ is the
restriction to $U_J$ of a section $f_{j_0}\colon U_{j_0} \to Z$,
and likewise $f_J(1)$ is the restriction to $U_J$ of
a section $f_{j_1}\colon U_{j_1} \to Z$.
\endpr

We shall also need the notion of a multiparameter homotopy of
$\cK(\cA)$-complexes. A suitable concept for this is the following.

%
%
\medskip\ni\bf 3.3 Definition.  \sl
Let $h\colon Z\to X$ and $\cA$ be as in def.\ 3.2, and let $k\in \Z_+$.
\item{(i)} A {\rm \holo\ $(\cK(\cA),k)$-prism}
(or a $k$-prism over $\cK(\cA)$) with values in $Z$ is a
collection 
$$  f_* = \{f_J \colon |J|\times [0,1]^k \to \cO_h(U_J,Z),
    \ \ J \in \cK(\cA) \},
$$
where $\cU$ is a faithful open \nbd\ of $\cA$, such that for
each fixed $y\in [0,1]^k$ the associated family
$$
   f_{*,y} = \{f_{J,y}=
   f_J(\cdotp,y) \colon |J| \to \cO_h(U_J,Z),\ \
   J \in \cK(\cA) \},
$$
is a \holo\ $\cK(\cA)$-complex. A {\rm continuous $(\cK(\cA),k)$-prism}
with values in $Z$ is a collection 
$f_*=\{ f_J \colon |J|\times [0,1]^k \to \cC_h(U_J,Z)\colon J \in \cK(\cA)\}$
such that $f_{*,y}$ is a continuous $\cK(\cA)$-complex for each
$y\in [0,1]^k$.

\item{(ii)}
A prism $f_*$ is {\rm sectionally constant} if there is an open
set $U \supset \cup_{i\ge 0} A_i$ such that the complex $f_{*,y}$
is represented by a section $f_y\colon U\to Z$ for each fixed
$y\in [0,1]^k$. If this holds for all $y$ in a subset $Y\subset [0,1]^k$,
we say that $f_*$ is sectionally constant on $Y$.
\medskip\rm

Thus a $(\cK(\cA),0)$-prism is a $\cK(\cA)$-complex, a
$(\cK(\cA),1)$-prism is the same thing as a homotopy of
$\cK(\cA)$-complexes, a 2-prism is a homotopy of 1-prisms,
etc. A sectionally constant $(\cK(\cA),k)$-prism is a $k$-parameter
homotopy of sections over an open \nbd\ of the union of sets in $\cA$.

%
%
%
%
\beginsection 4. Cartan strings and the initial \holo\ complex.

Let $X$ be a complex manifold. A compact set $C\subset X$,
contained in another compact set $B\subset X$, is said to be
{\it Runge} in $B$ if $C$ has a basis of open \nbd s each
of which is Runge in a fixed \nbd\ of $B$. A compact set
$K \subset X$ is said to be a {\it Stein compactum} if it
has a basis of open \nbd s which are Stein.

%
%
\proclaim 4.1 Definition:
An ordered pair of compact sets $(A,B)$ in a complex manifold
$X$ is said to be a {\rm Cartan pair} (or a Cartan string of
length 2) if
\item{(i)}    each of the sets $A$, $B$, and $A\cup B$ has a
basis of Stein \nbd s,
\item{(ii)}   $\bar{A\bs B} \cap \bar{B\bs A} =\emptyset$
\ \ (the separation condition), and
\item{(iii)}  the set $C=A\cap B$ is Runge in $B$.
($C$ may be empty.)

Our definition of a Cartan pair in [FP1] was similar, except that it
did not include the property (iii). We now recall Gromov's definition
of a {\it Cartan string}. It is by induction on $n\in N$,
a Cartan string of length $2$ being just a Cartan pair.

\medskip\ni\bf 4.2 Definition. \rm ([Gr1], 4.2.D'.) \sl
Let $X$ be a complex manifold and $A_0,A_1,\ldots,A_n \subset X$
compact subsets ($n\ge 1$).  The sequence  $(A_0,A_1,\ldots,A_n)$ 
is a {\rm Cartan string} of length $n+1$ if
\item{(i)}  $(A_0 \cup \ldots \cup A_{n-1}, A_n)$  is a Cartan
pair, and
\item{(ii)} if $n\ge 2$, the sequences $(A_0,\ldots,A_{n-1})$
and $(A_0 \cap A_n,\ldots,A_{n-1} \cap A_n)$ are Cartan strings
of length $n$.

\ni A locally finite covering $\cA=\{A_0,A_1,A_2,\ldots\}$ of $X$ by
compact sets is a {\rm Cartan covering} if $(A_0,A_1,\ldots,A_n)$
is a Cartan string for each $n\in \N$.
\medskip\rm

Note that each set $A_k$ in a Cartan string, and also each
finite union $A^{k}= \bigcup_{i=0}^{k} A_i$, is a Stein compactum 
in $X$. The order of sets in a Cartan string is important 
because of the property (iii) in Definition 4.1; hence a Cartan 
string is a {\it sequences} of sets and not merely a family. Cartan
strings have the following {\it hereditary property:}

%
%
\proclaim 4.3 Proposition: If $(A_0,A_1,\ldots,A_n)$ is a
Cartan string in a complex manifold $X$, and if $B\subset X$ is
a Stein compactum in $X$, then $(A_0\cap B,\ldots,A_n\cap B)$
is also a Cartan string.

\demo Proof: Induction on $n\in \N$. Note first that, if $A$ and
$B$ are compact sets in $X$ with bases of Stein \nbd s, then
$A\cap B$ also has a basis of Stein \nbd s (since intersections
of Stein domains are again Stein). Consider first the case $n=1$,
i.e., we have a Cartan pair $(A_0,A_1)$, and we wish to prove
that $(A_0\cap B,A_1\cap B)$ is also a Cartan pair.
Clearly it satisfies (i) and (ii) in def.\ 4.1.
Property (iii) follows from

\proclaim 4.4 Lemma: If $D_0\subset D_1$ and $\Omega$ are open
Stein domains in a complex manifold $X$, and if $D_0$ is
Runge in $D_1$, then $D_0\cap\Omega$ is Runge in $D_1\cap \Omega$.

\demo Proof of Lemma 4.4: Choose any compact set
$K\ss D_0\cap \Omega$. We denote by $\wh K_D$ the
\holo ally convex hull of $K$ with respect to a domain
$D\subset X$ containing $K$. Clearly
$\wh K_{D_1\cap \Omega} \subset \wh K_{D_1} \cap \wh K_\Omega$.
Since $D_0$ is Runge in $D_1$ and both domains are Stein,
we have $\wh K_{D_0}=\wh K_{D_1}$ ([H\"or], Theorem 2.7.3),
and therefore
$$
    \wh K_{D_1\cap \Omega} \subset \wh K_{D_1} \cap \wh K_\Omega
    = \wh K_{D_0} \cap \wh K_\Omega \ss D_0\cap \Omega.
$$
It follows (Theorem 2.7.3 in [H\"or]) that $D_0\cap \Omega$ is
Runge in $D_1\cap \Omega$.
\endpr

This completes the proof of Proposition 4.3 when $n=1$. Suppose
now that the result holds for some $n$. Let $(A_0,A_1,\ldots,A_{n+1})$
be a Cartan string of length $n+2$. To see that
$(A_0\cap B,\ldots,A_{n+1}\cap B)$ is a Cartan string,
we must verify that:
\item{(i)} the pair of sets
$ (A_0\cap B)\cup \cdots \cup (A_n \cap B)
   =(A_0\cup\cdots\cup A_n)\cap B
$
and $A_{n+1}\cap B$ is a Cartan pair. Since
$(A_0\cup\cdots\cup A_n,A_{n+1})$ is a Cartan pair,
this follows from the case $n=1$ proved above;
\item{(ii)} the strings $(A_0\cap B,\cdots,A_n\cap B)$ and
$(A_0\cap A_{n+1}\cap B,\ldots,A_n\cap A_{n+1}\cap B)$ are
Cartan strings of length $n+1$. This follows immediatelly
from Definition 4.3 and from the inductive hypothesis.
\endpr

%
%
%
%
\proclaim 4.5 Corollary: If $\cA=\{A_0,A_1,A_2,\ldots\}$ is
a sequence of compact sets in a complex manifold $X$ such that
for each $n\in \N$ the pair $(A_0\cup\cdots\cup A_{n-1},A_n)$ is
a Cartan pair, then  for each $n\in \N$ the string
$(A_0,A_1,\cdots,A_n)$ is a Cartan string.

\demo Proof: This follows from Proposition 4.3 by induction on $n$.
\endpr

%
%
%
%
\proclaim 4.6 Theorem: For each open covering $\cU=\{U_j\}$ of a
Stein manifold $X$ there exists a Cartan covering
$\cA=\{A_i \colon i=0,1,\ldots\}$ of $X$ which is
subordinate to $\cU$, i.e., such that each set $A_i$ is
contained in $U_j$ for some $j=j(i)$. Moreover, if $K\subset X$
is a compact \hc{X}\ subset and $U_0\subset X$ is an open set
containing $K$, we can choose $\cA$ such that
$K\subset A_0\subset U_0$ and $A_i\cap K=\emptyset$ for
$i\ge 1$.

\demo Proof: This was proved by Henkin and Leiterer in
sect.\ 2 of [HL]. We recall briefly the main idea.
The conditions imply that there exists a smooth \spsh\
exhaustion function $\rho\colon X\to \R$ with nondegenerate
critical points (a Morse function) such that
$K\subset \{\rho< 0 \} \ss U_0$ and $0$ is a regular value
of $\rho$. Set $A_0=\{\rho\le 0\}$. One can reach any higher
sublevel set $\{\rho\le c\}$, for $c>0$ being a regular value of
$\rho$, by successively adding (finitely many times) small
strongly pseudoconvex domains $A_k$ to the union of the previous
sets $A^{k-1}= \cup_{i=0}^{k-1} A_i$ in such a way that for
each $k$ the pair $(A^{k-1},A_k)$ is a {\it special \psc\ bump}
in the terminology of [HL]. It is clear from their definition
that such a pair is also a Cartan pair as defined here. Moreover,
we can insure that each set $A_k$ for $k\ge 1$ is contained in
a set $U_j$ for some $j\ge 1$, and the resulting family is
locally finite. Corollary 4.5 above implies that the
covering $\cA=\{A_0,A_1,A_2,\ldots\}$ of $X$ that one builds
in this way is a Cartan covering.

If $\rho$ has no critical values
on an interval $[c_0,c_1] \subset \R$ then we can actually reach
$\{\rho\le c_1\}$ from $\{\rho\le c_0\}$ by adding {\it convex bumps}
(this is called a noncritical pseudoconvex extension).
To cross the critical points of $\rho$ one must attach more
general pseudoconvex bumps.
\endpr

%
%
The following  proposition provides a homotopy of complexes
(a 1-prism) from the initial continuous section to a \holo\
complex over a Cartan covering of $X$.

\proclaim 4.7 Proposition:
{\rm (Construction of the initial holomorphic complex) }
Let $X$ be a Stein manifold, $K\subset X$ a compact \hc{X}\ subset
and $h\colon Z\to X$ a \holo\ submersion onto $X$ with the
property that each point $x\in X\bs K$ has an open \nbd\
$U_x\subset X$ such that $h\colon Z|U_x \to U_x$ admits a spray.
Let $a\colon X\to Z$ be a continuous section which is \holo\ in
an open set $U_0\supset K$. Then there exists a Cartan
covering $\cA=\{A_i \colon i=0,1,\ldots\}$ of $X$ and a
continuous $(\cK(\cA),1)$-prism $a_*=\{a_{*,s} \colon s\in [0,1]\}$
with values in $Z$ such that
\item{(i)}   $K\subset A_0 \subset U_0$, $K\cap A_i=\emptyset$
for $i\ge 1$, and $a_{(0),s}=a|U_0$ for each $s\in [0,1]$,
\item{(ii)}  $a_{*,0}$ is the constant $\cK(\cA)$-complex given by
the section $a \colon X\to Z$,
\item{(iii)} $a_{*,1}$ is a \holo\ $\cK(\cA)$-complex,
\item{(iv)}  for each $j\ge 1$ the submersion $Z\to X$ admits a
spray over an open set $U_j \supset A_j$.

\demo Proof:
Let $N=\dim Z = n+m$, where $n=\dim X$ and $m$ is the dimension of
the fibers $h^{-1}(x)$ ($x\in X$). Denote by $P^N$ the unit
open polydisc in $\C^N$ with complex coordinates $\z=(\z',\z'')$,
where $\z'\in \C^n$ and $\z''\in \C^m$. Let $\pi\colon P^N\to P^n$
be the projection $\pi(\z',\z'')=\z'$. Since $h\colon Z\to X$ is
a submersion, there exist for each point $z_0\in Z$ open \nbd s
$V \subset Z$ of $z_0$, $U\subset X$ of $x_0=h(z_0)$, and
biholomorphic maps $\Phi\colon V\to P^N$, $\phi\colon U\to P^n$,
such that $\pi\circ \Phi = \phi\circ h$ on $V$ and $\Phi(z_0)=0$.
Such map $\Phi$ induces a linear structure on the fibers of
$h|V$ which lets us add sections of $h\colon V\to U$ and
take their convex linear combinations.

If $z_0$ belongs to the graph of $a$, we can choose these
\nbd s and maps such that $a(U)\subset V$. In this case
$\Phi$ maps $a(U)$ onto the graph of a section
$\wt a(\zeta')=(\zeta',a''(\zeta'))$ ($\zeta'\in P^n$) of
the projection $\pi\colon P^N\to P^n$. The family
$a_s\colon U\to V$, given by
$$
    a_s(x)=\Phi^{-1}\bigl( \phi(x), (1-s)a''(\phi(x)) \bigr)
      \quad (x\in U,\ 0\le s\le 1),
$$
is a homotopy of continuous sections of $h$ over $U$ such
that $a_s(U)\subset V$ for each $s\in [0,1]$, $a_0=a|U$,
and the section $a_1$ is \holo. By shrinking $U$ around
$x_0$ and replacing $V$ by $V\cap h^{-1}(U)$ we may insure
in addition that the graph of the entire homotopy $a_s$ stays
in a prescribed open \nbd\ of $\bar{a(U)}$.

Let $U_0\subset X$ be an open set containing $K$ such that
$a$ is \holo\ in $U_0$. Set $V_0=h^{-1}(U_0) \subset Z$
and $a_{(0),s}=a|U_0$ for $s\in [0,1]$. Using the above argument
we can cover the graph of $a$ outside of $V_0$ by open
neighborhoods $V_j\subset Z$ ($j=1,2,3,\ldots$) biholomorphic
to $P^N$, with $U_j=h(V_j) \subset X$ biholomorphic to $P^n$,
such that for each $j\in \N$ we have a homotopy of
continuous sections $a_{(j),s} \colon U_j\to V_j$
($0\le s\le 1$) of $h$ satisfying
\item{(a)} $a_{(j),0}=a|U_j$,
\item{(b)} the section $a_{(j),1}$ is \holo\ in $U_j$,
\item{(c)} $U_j\cap K=\emptyset$ and $Z|U_j$ admits a spray
for each $j\ge 1$,
\item{(d)} if $U_{(i,j)}=U_i\cap U_j\ne\emptyset$ for some
$i,j\ge 0$  then $a_{(i),s}(U_{(i,j)}) \subset V_j$
for each $s\in [0,1]$.

The property (d) insures that for each pair of indices
$i,j\in\Z_+$ such that $U_{(i,j)}\ne\emptyset$ there is a
homotopy of sections $a_{(i,j),s}(t)\colon U_{(i,j)}\to Z$,
depending continuously on $t,s\in [0,1]$, such that
$a_{(i,j),s}(0)=a_{(i),s}|U_{(i,j)}$,
$a_{(i,j),s}(1)=a_{(j),s}|U_{(i,j)}$,
$a_{(i,j),0}(t)=a|U_{(i,j)}$, and the section
$a_{(i,j),1}(t)$ is holomorphic on $U_{(i,j)}$ for each
$t\in [0,1]$. We get $a_{(i,j),s}(t)$ by taking the convex
linear combinations in $t\in [0,1]$ of the sections
$a_{(i),s}$ and $a_{(j),s}$, restricted to $U_{(i,j)}$,
where the combinations are taken with respect to the linear
structure on the fibers of $h|V_i$ induced by $\Phi_i$
(or with respect to the linear structure on the fibers of
$h|V_j$ induced by $\Phi_j$). If one of the indices is zero,
say $i=0$ and $j>0$, we must use the linear structure on
$V_j$ induced by $\Phi_j$ since there is no such structure
on $V_0$.

Likewise, if $U_{(i,j,k)}\ne \emptyset$ for some
multiindex $J=(i,j,k)$, we can use the linear structure on
the fibers in any one of the sets $V_i$, $V_j$, $V_k$
to get a homotopy of sections $a_{J,s}(t) \colon U_J\to Z$,
with $t$ belonging to the standard 2-simplex
$\disc^2 \subset \R^2$, whose restriction to the sides
of the simplex equals the respective homotopy obtained
in the previous step. Continuing this way we build a
1-prism $a_*$ on the covering $\cU=\{U_0,U_1,U_2,\ldots\}$
of $X$.

By Theorem 4.6 there is a Cartan covering
$\cA=\{A_0,A_1,A_2,\ldots\}$ of $X$ subordinate to $\cU$,
with $K\subset A_0\subset U_0$. Then $a_*$ induces in a
natural way a $(\cK(\cA),1)$-prism with the required properties.
This proves Proposition 4.7.
\endpr

%
%
%
%
\beginsection 5. Modifying holomorphic prisms over Cartan strings.

This section contains the heart of the proof of Theorems 1.2
and 1.5. In Proposition 5.1 we
show by induction on $n$ how to modify a \holo\
$\cK(\cA)$-complex $f_*$ over any finite subcomplex
$\cK^n=\cK(A_0,\cdots,A_n)$ into a \holo\ section, defined over
a \nbd\ of $A^n=A_0\cup A_1\cup\cdots\cup A_n$, provided
that $(A_0,\ldots,A_n)$ is a Cartan string. Since the inductive
procedure requires us to solve the problem for parametrized
families of complexes, we consider \holo\ prisms from the
outset. The initial case is $n=1$ when the string
$(A_0,A_1)$ is a {\it Cartan pair}; here we need the
analytic tools from [FP1].

%
%
%
%
\medskip\ni \bf 5.1 Proposition. \sl
Let $h\colon Z \to X$ be a \holo\ submersion. Let $(A_0,\ldots,A_n)$
be a Cartan string in $X$ and $\cK^n= \cK(A_0,\ldots,A_n)$ its nerve.
Assume that for each $i=1,\ldots n$ there is an open set $U_i \supset A_i$
in $X$ such that $h\colon Z|U_i \to U_i$ admits a spray.  If $f_*$ is a
holomorphic $(\cK^n,k)$-prism with values in $Z$ which is
sectionally constant on a nice compact set $Y\subset [0,1]^k$
(def.\ 1.4), there is a homotopy $f^u_*$ $(0\le u\le 1)$ of
\holo\ $(\cK^n,k)$-prisms such that
\item{(i)}   $f^0_* = f_*$ is the given prism,
\item{(ii)}  the prism $f^1_*$ is sectionally constant,
\item{(iii)} for each $y\in [0,1]^k$ and $u\in [0,1]$,
the section $f^u_{(0),y}$ approximates $f_{(0),y}$ on $A_0$
as well as desired,
\item{(iv)}  $f^u_{*,y} = f_{*,y}$ for all $y \in Y$ and
$u\in [0,1]$ (i.e., the homotopy is fixed on $Y$).

\ni Moreover, if the restriction $f_*|\cK^{n-1}$ to the subcomplex
$\cK^{n-1}=\cK(A_0,\ldots,A_{n-1})$ is sectionally constant,
the homotopy $f^u_*$ can be chosen such that, in addition to
the above, the prism $f^u_*|\cK^{n-1}$ is sectionally constant
for each $u\in [0,1]$ and the corresponding sections
$f^u_{*,y}|\cK^{n-1}$ for $y\in [0,1]^k$, which are \holo\
in a \nbd\ of $A^{n-1}$, approximate $f_{*,y}|\cK^{n-1}$
uniformly on $A^{n-1}$.

\demo Proof: Replacing $X$ by a suitable Stein \nbd\
of $A^n$ we may assume it is Stein. The proof is by induction
on $n\ge 0$, and for $n=0$ there is nothing to prove.

\ni\it The case $n=1$.  \rm  Our data consists of a Cartan
pair $(A_0,A_1)$ in $X$ and a \holo\ $(\cK(A_0,A_1),k)$-prism
$f_*$ which is sectionally constant on a nice compact set
$Y\subset [0,1]^k$. Such an object $f_*$ is determined by
the following data:

\item{(a)} a pair of open sets $U_0 \supset A_0$ and
$U_1\supset A_1$,
\item{(b)} families of \holo\ sections
$a_y= f_{(0),y} \colon U_0\to Z$, $b_y=f_{(1),y} \colon U_1\to Z$,
depending continuously on $y\in [0,1]^k$,
\item{(c)} a family of \holo\ sections
$$
   c_{y,t}=f_{(0,1),y}(t)\colon U_{(0,1)}=U_0\cap U_1\to Z
$$
depending continuously on $t\in [0,1]$ and $y\in [0,1]^k$,

\ni such that $a_y|U_{(0,1)}=c_{y,0}$,\ \
$b_y|U_{(0,1)}=c_{y,1}$, and for each $y\in Y$
the section $c_{y,t}$ is independent of $t\in [0,1]$.
Hence for $y\in Y$ the family
$\{c_{y,t}\colon t\in[0,1]\}$ determines a \holo\
section $c_y\colon U_0\cup U_1\to Z$ such that $c_y|U_0 =a_y$
and $c_y|U_1=b_y$. We shall write $f_*=(a_*,b_*,c_*)$,
where $*$ indicates the missing parameters.

Our goal is to construct a homotopy $f^u_*=(a^u_*,b^u_*,c^u_*)$
($0\le u\le 1$) of \holo\ $(\cK(A_0,A_1),k)$-prisms
over smaller sets $U_0 \supset A_0$ and $U_1\supset A_1$
such that $f^0_*=f_*$ and $f^1_*$ is a constant prism,
which really means that $f^1_*$ is a collection
of \holo\ sections $f^1_y\colon U_0 \cup U_1\to Z$
(we have eliminated the $t$ parameter!). Moreover, the homotopy
must be fixed for $y\in Y$ and it must approximate the
sections $a_y$ over $A_0$ for each $y \in [0,1]^k$.
We shall denote the data in the homotopy $f^u_*$ by the same
letters as above, adding the upper index $u$.

The homotopy $f^u_*$ will be constructed in two steps.
For convenience we use the parameter interval $u\in [0,2]$
and later rescale it to $[0,1]$. In the first step we apply
the h-Runge theorem (Theorem 4.2 in [FP1]) to get a homotopy
$f^u_*$, $0\le u\le 1$, from $f^0_*=f_*$ to another prism
$f^1_*$ such that we do not move the section $a_y$
(i.e., $a^u_y=a_y\colon U_0\to Z$ for all $u$ and $y$),
and such that the section $b^1_y\colon U_1\to Z$ approximates
$a_y$ in a \nbd\ of $A_0\cap A_1$ for each $y\in [0,1]^k$.
In the second step we apply the gluing lemma, Theorem 5.5 in [FP1],
to obtain homotopies of sections
$$
    a^u_y \colon U_0\to Z, \quad b^u_y \colon U_1\to Z,\quad
    c^u_{y,t}\colon U_0\cap U_1\to Z \qquad (1\le u\le 2)
$$
such that at $u=2$, $a^2_y=b^2_y$ over $U_0\cap U_1$ for
each $y\in [0,1]^k$; hence these two sections define a single
\holo\ section $f^2_y\colon U_0 \cup U_1\to Z$.
Of course we shrink the sets $U_0$ and $U_1$ again when gluing.
Moreover, both homotopies will be fixed on $Y$.

Consider the first step. Since $A_0\cap A_1$ is Runge in $A_1$
and the submersion $h\colon Z\to X$ admits a spray over a
\nbd\ of $A_1$, the h-Runge theorem with parameters (Theorem 4.2
in [FP1] and the remark following it) shows that, after shrinking
the sets $U_0$ and $U_1$, there exists a homotopy of \holo\ sections
$g^s_{y,t}\colon U_0\cap U_1 \to Z$ ($0\le s\le 1$), depending
continuously on $t,s,y$, satisfying
\item{(i)}   $g^0_{y,t}=c_{y,t}$ for each $y$ and $t$,
\item{(ii)}  $g^s_{y,1}=c_{y,1}=b_y|U_0\cap U_1$ for each $s$
and $y$,
\item{(iii)} $g^1_{y,t}$ extends to a \holo\ section over $U_1$
for each $y$ and $t$,
\item{(iv)}  the homotopy is fixed on $Y$, i.e., for $y\in Y$ we have
$g^s_{y,t}=c_y|U_0\cap U_1$ for each $s$ and $t$, and
\item{(v)}   $g^s_{y,t}$ approximates $c_{y,t}$ in a \nbd\ of $A_0\cap A_1$
as close as desired, uniformly with respect to all parameters.

\medskip\ni
We define $f^u_*=(a^u_*,b^u_*,c^u_*)$ for $0\le u\le 1$ by
$$  \eqalign{
    a^u_y     &= a_y, \quad  b^u_y = g^1_{y,1-u}, \cr
    c^u_{y,t} &= \cases{ c_{y,t}    & if $0\le t\le 1-u$; \cr
                         {\phantom{----}} & \cr
                g^{(t+u-1)/u}_{y,1-u} & if $1-u < t\le 1$. \cr}
         \cr}
$$
The reader will easily verify that this satisfies all required
properties. In particular, at $u=1$, $b^1_y=g^1_{y,0}$ approximates
$c_{y,0}=a_y$ in a \nbd\ of $A_0\cap A_1$.

Next we apply Theorem 5.5 in [FP1] (on gluing parametrized
families of holomorphic sections over Cartan pairs) to get
homotopies of sections $a^u_y \colon U_0\to Z$ and
$b^u_y \colon U_1\to Z$ for $1\le u\le 2$ such that $a^u_y$
approximates $a^1_y=a_y$ on $A_0$ for each $u\in [1,2]$
and $a^2_y=b^2_y$ in $U_0\cap U_1$. Moreover, over a \nbd\ of the
set $A_0$, the graphs of all sections $a^u_y$, $b^u_y$
($1\le u\le 2$) and $c^1_{y,t}$ ($0\le t\le 1$) lie in small
\nbd\ of the (graph of) the section $a_y$ in $Z$, and such a
\nbd\ is biholomorphic to a \nbd\ of the zero section in a \hvb\
over $U_0$ (Lemma 5.3 in [FP1]). Using the resulting vector
bundle structure on this \nbd, we see that the triangle
of homotopies formed by these families is contractible,
meaning that it can be filled by a 2-parameter homotopy
$c^u_{y,t}$ ($0\le t\le 1$, $1\le u\le 2$) over a \nbd\
of $A_0\cap A_1$. This completes the proof of
Proposition 5.1 for $n=1$.

%
%
\medskip
\ni\it The induction step $n\Rightarrow n+1$. \rm
Suppose that Proposition 5.1 holds for all Cartan strings of
length $n+1$ for some $n\ge 1$ and for all $k\ge 0$. Let
$\cA=(A_0,\ldots,A_{n+1})$ be a Cartan string of length $n+2$
with the nerve $\cK^{n+1}=\cK(\cA)$, and  let $f_*$ be a
\holo\ $(\cK^{n+1},k)$-prism with values in $Z$ which is
sectionally constant on a nice compact subset
$Y \subset [0,1]^k$. Let
$\cK^n= \cK(A_0,\ldots,A_n) \subset \cK^{n+1}$. The proof
consists of the following three steps, each of which is
accomplished by constructing a suitable homotopy of prisms.

\medskip
\ni \it Step 1: \rm
Reduction to the case when $f_*|\cK^n$ is a sectionally
constant prism.

\ni\it  Step 2: \rm Reduction to the case when $f_*$ represents
a $(k+1)$-prism over the Cartan pair $(A^n,A_{n+1})$, where
$A^n=A_0\cup A_1\cup\cdots \cup A_n$.

\ni \it Step 3: \rm Applying the case $n=1$ to the prism in step 2
to get a sectionally constant $(\cK^{n+1},k)$-prism.

\medskip
We begin by some general considerations. We denote the coordinates
on $\R^{n+1}$ by $t=(t',t_{n+1})$, where $t'=(t_1,\ldots,t_n)\in \R^n$,
and identify $\R^n$ with the coordinate hyperplane
$\R^n\times \{0\} \subset \R^{n+1}$. The body $K^{n+1}$
of the nerve $\cK^{n+1}$ can be represented as the union of
certain faces of the standard simplex $\disc^{n+1}\subset \R^{n+1}$.
(In fact $K^{n+1}=\disc^{n+1}$ \iff\
$A_0\cap A_1\cap\cdots\cap A_{n+1}\ne \emptyset$.)
The body $K^n \subset \R^n$ of the subcomplex
$\cK^n=\cK(A_0,\cdots,A_n)\subset\cK^{n+1}$ is precisely the
base $K^{n+1} \cap \{t_{n+1}=0\}$ of $K^{n+1}$. We shall also
need the complex
$$
    \cK^n_1 = \cK(A_0 \cap A_{n+1}, \ldots, A_n \cap A_{n+1})
      \subset \cK^n.  \eqno(5.1)
$$
Note that $\cK^n_1=\{J\in \cK^n\colon (J,n+1)\in \cK^{n+1}\}$.
Its body $K^n_1$ is a subset of $K^n$ which equals
$\bar{(K^{n+1}\bs K^n)} \cap K^n$.
Moreover, for each $0<s<1$, the section
$K^{n+1}\cap \{t_{n+1}=s\}$ is homeomorphic
to $K^n_1$. The map
$$
    r\colon \R^{n+1}\to\R^{n+1},\quad r(t,s)=(t(1-s),s)\quad
         (t\in \R^n,\ s\in \R)  
$$
maps the prism $\disc^n\times [0,1] \subset \R^{n+1}$ onto
the standard simplex $\disc^{n+1}$ (it is homeomorphic outside
the level $s=1$), and it maps $K^n_1 \times \{s\}$
homeomorphically onto the section $K^{n+1}\cap \{t_{n+1}=s\}$
for each $s\in (0,1)$.

\medskip\ni \it Proof of step 1. \rm
Since $\wt f^0_* = f_*|\cK^n = \{ f_J \colon J \in \cK^n \}$
is a $k$-prism over a Cartan string of length $n+1$, the
induction hypothesis provides a \holo\ homotopy
$\wt f_*= \{\wt f^{u}_* \colon u\in [-1,0]\}$
such that each $\wt f^u_*$ is a $(\cK^n,k)$-prism, the homotopy
is fixed for all $y\in Y$, and the prism $\wt f^{-1}_*$
is sectionally constant. The parameter space of the prism $f_*$
is $K^{n+1}\times [0,1]^k$ while the parameter space of
$\wt f_*$ is $K'\times [0,1]^k$, where
$$
    K'= \{(t',u) \in \R^n\times\R \colon\
             t'\in K^n,\ -1\le u\le 0 \}.
$$
Note that $f_*$ and $\wt f_*$ agree on the intersections
of their domains $K^n \times [0,1]^k$ and hence define a collection
of sections, parametrized by the set $L\times [0,1]^k$, where
$L=K^{n+1}\cup K' \subset \R^{n+1}$. We denote this collection
by $\{g_y(t) \colon t\in L,\ y\in [0,1]^k \}$.

For each $s\in [0,1]$ we denote by $L_s \subset \R^{n+1}$ the set
$$
    L_s= (K^{n+1}\bs K^n) \cup
    \{ (t',t_{n+1}) \colon t'\in K^n_1,\ -s\le t_{n+1} \le 0\}
    \cup \{(t',-s)\colon t'\in K^n\} .
$$
Intuitively speaking, $L_s$ is obtained by pushing the base $K^n$ of
$K^{n+1}$ for $s$ units in the negative $t_{n+1}$ direction and
then adding to this the vertical sides $K^n_1\times [-s,0]$.
Clearly $L_0=K^{n+1}$, and $L_s$ is homeomorphic to $K^{n+1}$
for each $s\in [0,1]$. In fact, there is a continuous family
of homeomorphisms
$\Theta_s\colon K^{n+1}\to L_s$ ($0\le s\le 1$) such that
$\Theta_0$ is the identity, each $\Theta_s$ preserves the top
vertex $(0,\ldots,0,1)\in K^{n+1}$ and the cellular structure of
the two sets, and $\Theta_s$ maps $K^n$ (the base of $K^{n+1}$)
onto $K^n\times \{-s\}$ (the base of $L_s$) by a downward shift
for $s$ units. By `respecting the cellular structure' we mean
the following. Each face $J\in \cK^n_1$ determines a face
$\wt J=(J,n+1)\in \cK^{n+1}$, and $\Theta_s$ maps its body
$|\wt J| \subset K^{n+1}$ onto the set
$|\wt J|\cup \{ (t',t_{n+1}) \colon
  t'\in |J|,\ -s\le t_{n+1}\le 0 \} \subset L_s$.
Using this family we define a homotopy $H^u_*$ ($0\le u\le 1$)
of $(\cK^{n+1},k)$-prisms by
$$
    H^u_J(t) =  g_y(\Theta_{u}(t)) \quad
    (J\in \cK^{n+1},\ t\in |J| \subset K^{n+1}).
$$
Clearly $H^0_*=f_*$ and $H^1_*|\cK^n = \wt f^{-1}_*|\cK^n$
is sectionally constant. This proves step 1.

\medskip \ni \it
Proof of step 2: \rm
By step 1 we may assume that the prism $f_*$ is such that
$f_*|\cK^n$ is sectionally constant. The next step is to modify
$f_*$ by a homotopy of \holo\ prisms into another prism which
is sectionally constant also in the direction of the last
variable $t_{n+1}$. Let $\cK^n_1$ be the complex (5.1).
We associate to $f_*$ a \holo\ $(\cK^n_1,k+1)$-prism
$$
    F_* =\{F_{J,(y,s)}\colon |J|\to \cO_h(U_{(J,n+1)},Z),\quad
       J\in \cK^n_1,\ y\in [0,1]^k,\ s\in [0,1] \}
$$
where
$$
    F_{J,(y,s)}(t) = f_{(J,n+1),y}( r(t,s)) \quad
            (t \in |J|,\ y\in [0,1]^k,\ s\in [0,1]).
$$
The set
$$
    Y_1 = (Y \times [0,1]) \cup ([0,1]^k \times \{0,1\})
             \subset [0,1]^{k+1}
$$
is a nice compact subset of $[0,1]^{k+1}$ (def.\ 1.4).
Since $f_*|\cK^n$ is sectionally constant, the prism $F_*$
(which is associated to the complex $\cK^n_1$ of length $n+1$)
satisfies the induction hypothesis with respect to the set
$Y_1$. Hence there is a homotopy $F^u_*$ ($u \in [0,1]$) of
\holo\ $(\cK^n_1,k+1)$-prisms, beginning with $F^0_*=F_*$,
such that the homotopy is fixed for $(y,s)\in Y_1$ and such
that $F^1_*$ is sectionally constant. This means that for each
fixed $(y,s) \in [0,1]^{k+1}$ the $\cK^n_1$-complex
$F^1_{*,(y,s)}$ is constant, i.e., it represents a \holo\
section $F^1_y(s)\colon V\to Z$ of $Z\to X$
over an open set $V \supset A^n\cap A_{n+1}$.
In particular, since the homotopy is fixed for $s=0$ and $s=1$,
the section $F^1_y(0)$ coincides with the section represented by
the constant complex $f_{*,y}|\cK^n$, and the section $F^1_y(1)$
coincides with the section $f_{(0,\ldots,0,1),y}$
associated to the set $A_{n+1}$.

\medskip\ni \it
Proof of step 3: \rm
We may consider the family of sections
$$
    F^1_*=\{F^1_y(s) \colon V\to Z |\ s\in [0,1],\ y\in [0,1]^k\},
$$
obtained in step 2, as a \holo\ $k$-prism over the complex
$\cK'=\cK(A^n,A_{n+1})$ determined by the Cartan pair
$(A^n,A_{n+1})$. Here the parameter $s\in [0,1]$ represents
the variable in the body $|\cK'|=[0,1]$. For each $y\in [0,1]^k$
the section $F^1_y(0)$ extends holomorphically to a \nbd\ of $A^n$
and $F^1_y(1)$ extends to a \nbd\ of $A_{n+1}$.

The case $n=1$ of Proposition 5.1 gives us a homotopy
$G^u_*$ ($0\le u\le 1$) of \holo\ $(\cK',k)$-prisms
such that $G^0_*=F^1_*$, $G^1_y$ is a constant $\cK'$-complex
for each $y\in [0,1]^k$ (i.e., a \holo\ section over an
open \nbd\ of $A^{n+1}=A_0\cup\cdots\cup A_{n+1}$),
and the homotopy is fixed for $y\in Y$ (where $G^0_y=F^1_y$
is already a section over $A^{n+1}$). Moreover, for each
$u\in [0,1]$ the section $G^u_y(0)$ (which is \holo\ over a \nbd\
of $A^n$) approximates the section $F^1_y(0)=f_{*,y}|\cK^n$ on
the set $A^n$.

We can interpret the families $\{F^u_* \colon 0\le u\le 1\}$
and $\{G^u_* \colon 0\le u\le 1\}$ as homotopies of
\holo\ $(\cK^{n+1},k)$-prisms. By connecting these homotopies
$F^u_*$ and $G^u_*$ (in this order) we obtain a homotopy $f^u_*$
($0\le u\le 1$) of \holo\ $(\cK^{n+1},k)$-prisms, beginning at
$u=0$ with $f_*$ and ending at $u=1$ with the sectionally
constant prism $G^1_*$.

If we assume in addition that the restriction $f_*|\cK^n$
is sectionally constant on $[0,1]^k$ (so $f_{*,y}|\cK^n$ is
a \holo\ section in a \nbd\ of $A^n$ for each $y\in [0,1]^k$),
we can skip the initial step in the proof of the inductive step.
Note that, by construction, the restriction $F^u_*|\cK^n$ is
independent of $u \in [0,1]$ (since the homotopy $F^u_*$ is fixed
on the set $Y_1$), and the homotopy $G^u_*$ is such that the
complex $G^u_{*,y}|\cK^n$ is represented by a \holo\ section
in a \nbd\ of $A^n$ which approximates the section
$F^1_{*,y}|\cK^n = f_{*,y}|\cK^n$ on $A^n$, uniformly with respect
to $u\in [0,1]$ and $y\in [0,1]^k$. Hence the section
$f^u_{*,y}|\cK^n$ approximates $f_{*,y}|\cK^n$ on $A^n$,
the approximation being uniform with respect to $u\in [0,1]$
and $y\in [0,1]^k$. This concludes the proof of the 
induction step.
\endpr

The next proposition shows that a \holo\ 1-prism can be extended
from a finite subcomplex to the entire complex such that the
$0$-level of the prism matches a given complex. This does not
require any analytic tools and hence the result applies to
any locally finite family $\cA$.
%
%
%
%
\proclaim 5.2 Proposition: Let $\cA=\{A_0,A_1,A_2,\ldots\}$ be
a locally finite family of compact sets in a complex manifold
$X$. Denote its nerve by $\cK(\cA)$ and let
$\cK^n=\cK(A_0,\ldots,A_n) \subset \cK(\cA)$ for each $n\in \N$.
Assume that $h\colon Z\to X$ is a \holo\ submersion onto $X$,
$f_*$ is a \holo\ $\cK(\cA)$-complex with values in $Z$,
and $g_*$ is a \holo\ $(\cK^n,1)$-prism for some $n\in \N$
such that $g_{*,0}=f_*$. Then there exists a \holo\
$(\cK(\cA),1)$-prism $G_*$ such that
$G_{*,0}=f_*$ and $G_*|\cK^n =g_*$.

\demo Remark: A similar result holds if $f_*$ is a \holo\
$(\cK(\cA),k)$-prism and $g_*$ is a \holo\ $(\cK^n,k+1)$-prism
with base $f_*|\cK^n$: such $g_*$ extends to a \holo\
$(\cK(\cA),k+1)$-prism $G_*$.

\demo Proof:
We choose representatives of $f_*$ and $g_*$
defined on a faithful open \nbd\ $\cU$ of $\cA$ (see def.\ 3.2).
Write $A^n=A_0\cup\cdots\cup A_n$ as before. Let $m\ge n$ be the
smallest integer such that $A_k \cap A^n = \emptyset$ for each
$k\ge m$. We represent the set $K^m=K(A_0,\ldots,A_m)$
(the body of the subcomplex $\cK^m \subset \cK(\cA)$)
as a subset of $\R^m$. We denote the coordinates on $\R^{m+1}$
by $(t,s)$, with $t\in \R^m$ and $s\in \R$, and identify $\R^m$
with $\R^m\times\{0\} =\{s=0\} \subset \R^{m+1}$.
Similarly we identify a set $K\subset \R^m$ with
$K\times \{0\}\subset \R^{m+1}$ and write
$K\times [0,1]=\{(t,s)\colon t\in K,\ s\in [0,1]\}$.
For each face $J\in \cK^m$ we denote by $b|J|\subset K^m$
the boundary of its body $|J|$.

%
%
\proclaim 5.3 Lemma: There exists a retraction
$$
    r\colon K^m\times [0,1]  \to
      K^m \cup (K^n\times [0,1]) \subset \R^{m+1}
$$
such that for each face $J\in \cK^m\bs \cK^n$ we have
\item{(i)}  $r(|J|\times [0,1])\subset |J|\cup (b|J|    \times [0,1])$,
\item{(ii)} if $|J|\cap K^n=\emptyset$ then $r(t,s)=t$ for each
$t\in |J|$ and $s\in [0,1]$.

\demo Proof: We first define $r$ over those faces $J \in \cK^m$
for which either $|J|\subset K^n$ (we let $r$ be the identity
on $|J|\times [0,1]$) or $|J|\cap K^n=\emptyset$ (we let $r(t,s)=t$
for $t\in |J|$). We also define $r$ to be the identity map on
the bottom side $K^m=K^m\times\{0\}$. On the remaining faces
$J\in \cK^m$ we define $r$ inductively with respect to the
dimension of $J$. Suppose that $r$ has already been defined on
all faces of dimension $<k$ and let
$J=(j_0,\ldots,j_k) \in \cK^m$. Then $r$ is defined on
$|J|\cup (b|J|\times [0,1])$ and satisfies (i), and it satisfies
(ii) on those sides of $b|J|$ which are disjoint from $K^n$.
Moreover, $r$ is the identity on $|J|=|J|\times\{0\}$.
It is now clear that $r$ extends from
$|J|\cup (b|J|\times [0,1])$ to $|J|\times [0,1]$ so that (i) holds.
\endpr

Let $r$ be as in Lemma 5.3. Write $r(t,s)=(r_0(t,s),u(t,s))$, where
$r_0(t,s) \in K^m$ and $u(t,s) \in [0,1]$. We define a \holo\
$(\cK^m,1)$-prism $G_*$ by setting for each $J\in \cK^m$, $t\in |J|$
and $s\in [0,1]$
$$
    G_{J,s}(t) = \cases{ f_{J}(r_0(t,s))  & if $u(t,s)=0$; \cr
                  g_{J,u(t,s)}(r_0(t,s)) & if $u(t,s)>0$. \cr }
$$
Property (i) in Lemma 5.3 implies that the section $G_{J,s}(t)$
for $t\in |J|$ is defined (and holomorphic) in the set $U_J$
(it may be holomorphic in a larger set if $r_0(t,s)\in b|J|$,
but in such case we restrict it to $U_J$). The collection
$G_*=\{G_{J,s}\colon J\in \cK^m,\ s\in [0,1]\}$
is then a \holo\ $(\cK^m,1)$-prism which extends $g_*$ and
satisfies $G_{*,0}=f_*$.

The property (ii) of the retraction $r$ lets us extend $G_*$ to
a prism over the entire complex $\cK(\cA)$ by observing that
for those faces  $J\in \cK(\cA)$ which do not belong to $\cK^m$
we have $|J|\cap K^n=\emptyset$ (by definition of $m$) and
therefore $r(t,s)=t$ for $t\in |J|\cap K^m$. Thus we can simply
take $G_{J,s}(t)=f_{J}(t)$ for $t\in |J|$ and $s\in [0,1]$.
This completes the proof of Proposition 5.2.
\endpr

%
%
%
%
\beginsection 6. Proof of Theorem 1.5.

In this section we prove Theorem 1.5, using the tools developed
in section 5. We shall concentrate on the case of a single section;
the proof in the general parametric case is essentially the same.
Thus, we are given a continuous section $a\colon X\to Z$ which
is \holo\ in an open set $U_0 \supset K$, where $K\subset X$
is \hc{X}; our goal is to construct a homotopy
$H_s\colon X\to Z$ $(0\le s\le 1)$ of continuous sections such
that $H_0=a$, the section $H_1$ is holomorphic on $X$, and for each
$s\in [0,1]$ the section $H_s$ is \holo\ near $K$ and it approximates
$a$ on $K$.

Let $\cA=\{A_0,A_1,\ldots\}$ be a Cartan covering of $X$ given by
Theorem 4.6 such that $K\subset A_0\subset U_0$ and
$K\cap A_i=\emptyset$ for $i\ge 1$.
Also let $a_*=\{a_{*,s} \colon 0\le s\le 1 \}$ be a continuous
$(\cK(\cA),1)$-prism provided by Proposition 4.7. Thus the complex
$a_{*,0}$ is constant and represents the section $a$, $a_{*,1}$ is a
{\it \holo} $\cK(\cA)$-complex, and $a_{(0),s}=a|U_0$ for each
$s\in [0,1]$.

Let $d$ be a complete metric on $Z$ compatible with the
manifold topology. Fix an $\e>0$. We inductively construct
a sequence of \holo\ $\cK(\cA)$-complexes $f^n_*$ and
\holo\ $(\cK(\cA),1)$-prisms
$G^n_*=\{G^n_{*,s}\colon 0\le s\le 1\}$ for $n=0,1,2,\ldots$,
satisfying the following properties:

\medskip
\item{(a)} $f^0_* = a_{*,1}$,

\item{(b)} $G^n_{*,0}=f^n_*$ and $G^n_{*,1}=f^{n+1}_*$ for each
$n\in \Z_+$,

\item{(c)} for each $k\in \Z_+$, $n\ge k$ and $s\in [0,1]$
the complexes $f^n_*|\cK^k$ and $G^n_{*,s}|\cK^k$ are constant,
i.e., they are given by \holo\ sections denoted by
$f^n$ resp.\ $G^n_s$ in an open \nbd\ of
$A^k=A_0\cup\cdots\cup A_k$,

\item{(d)} (approximation) for each $n\in \Z_+$ and $s\in [0,1]$
we have
$$
    d\bigl( G^n_s(x), f^n(x) \bigr) < \e/2^{n+1} \quad (x\in A^n).
$$
In particular we have $d\bigl( f^{n+1}(x), f^n(x) \bigr)< \e/2^{n+1}$
for $x\in A^n$.

\medskip
In (d) we are using the notation for sections established in (c).
The property (d) implies that the sequence of sections
$f^n \colon A^n\to Z$ ($n=0,1,2,\ldots$) converges uniformly
on compacts in $X$ to a \holo\ section
$f^\infty=\lim_{n\to\infty} f^n \colon X\to Z$ which satisfies
$$
    d\bigl( f^\infty(x), a(x) \bigr) =
    d\bigl( f^\infty(x), f^0(x) \bigr) < \e
    \quad (x\in A_0).
$$
To construct a homotopy $H_s \colon X\to Z$ ($0\le s\le 1$)
between $H_0=a$ and $H_1=f^\infty$ we first construct a continuous
$(\cK(\cA),1)$-prism $h_*$ such that $h_{*,0}=a$ and
$h_{*,1}=f^\infty$. To do this we simply collect all individual
1-prisms $a_*$ and $G^n_*$ ($n\in \Z_+$) into a single 1-prism
as follows. For each $n\in \Z_+$ set $I_n=[1-2^{-n}, 1-2^{-n-1}]$
and let $\lambda_n\colon I_n\to [0,1]$ be the linear bijection
$\lambda_n(s)= 2^{n+1}(s-1+2^{-n})$. Then
$\cup_{n=0}^\infty I_n = [0,1)$. For $s\in [0,1)$ we define
$$
    h_{*,s}=\cases{ a_{*,2s}   & if $s\in I_0=[0,1/2]$; \cr
            G^{n-1}_{*,\l_n(s)} & if $s\in I_n,\ n\ge 1$.\cr}
$$
The two definitions of $h_{*,s}$ at the values $s=1-2^{-n}$
($n\in \Z_+$) are compatible  by (b). Properties (c)
and (d) imply that $\lim_{s\to 1} h_{*,s} = f^\infty$ uniformly
on compacts in $X$. Indeed, each compact set $L\ss X$ is contained
in some $A^m$, and for $n\ge m$ the complex $G^n_{*,s}$ is
constant on $A^n$, i.e., represented by a \holo\ section. Hence for
$1-2^{-n-1} \le s <1$ the complex $h_{*,s}|\cK^n$ is a \holo\
section in a \nbd\ of $A^n$. As $s\to 1$, these sections
converge uniformly on $A^n$ (and hence on $L$) to $f^\infty$.
This proves that, if we set $h_{*,1}=f^\infty$,
the collection $h_{*}=\{ h_{*,s} \colon 0\le s\le 1\}$
is indeed a continuous $(\cK(\cA),1)$-prism. Notice also that
the restriction of $h_*$ to $A_0$ (more precisely, to the
complex $\cK(A_0)$ represented by the first set $A_0$) is in
fact a homotopy of holomorphic sections $h_s$ ($0\le s\le 1$)
in a \nbd\ of $A_0$, connecting $h_0=a$ to $h_1=f^\infty$
and such that all sections in the family satisfy
$d(h_s(x),a(x))<\e$ for $x\in A_0$.

To complete the proof of Theorem 1.5 we apply the version
of Proposition 5.1 for continuous prisms to modify the 1-prism
$h_*$ by a homotopy of 1-prisms (keeping the ends $s=0$ and
$s=1$ fixed) into a 1-prism $H_*$ which is {\it sectionally
constant}, i.e., such that $H_*$ represents a homotopy of
continuous sections $\{H_s \colon X\to Z\colon 0\le s\le 1\}$.
Moreover we can achieve that in a \nbd\ of $A_0$ the two
sections $H_s$ and $h_s$ agree.

This concludes the proof of Theorem 1.5 in the case without
parameters. The parametric case is proved by the same tools
by introducing the parameter space $P$ into the definition of
\holo\ (and continuous) complexes and prisms and repeating
the same arguments in this setting. The analytic tools
used in the proof, namely the h-Runge theorem and the gluing
theorem, have been established in this generality in [FP1].

\demo A final remark:  We have described a procedure
in which the final holomorphic section of $Z\to X$ is obtained
as a locally uniform limit of holomorphic sections, defined
over increasingly larger compacts in $X$. Often one need less:
To modify a given continuous section, which is holomorphic
over a Stein compactum $K\subset X$, into another section that
is holomorphic over a larger Stein compactum $L\supset K$ in $X$.
Our proof accomplishes this in a {\it finite number of steps},
provided that {\it $L$ is a strongly pseudoconvex
extension of $K$}. This means that $K$ and $L$ are two
regular sublevel sets of a smooth \spsh\ function
defined in a \nbd\ of $\bar {L\bs K}$. In such case the
construction in [HL] gives a finite Cartan string
$(A_0,A_1,\ldots,A_n)$ in $X$ with $A_0=K$ and
$A^n= \bigcup_{0\le j\le n} A_j=L$. It suffices to
apply our proof on this finite string.

%
%
%
%
\beginsection 7. Proof of Theorem 1.6.

\demo Proof of Theorem 1.6 (a):
For each nonzero vector $v\in\C^q$ we denote by $[v] \in\CP^{q-1}$
the complex line in $\C^q$ spanned by $v$, and we denote
by $\pi_v\colon \C^q\to v^\perp = \C^{q-1}$ the orthogonal
projection onto $v^\perp$ with $\pi(v)=0$.

%
%
\proclaim 7.1 Lemma: {\rm (Existence of sprays.) }
Let $U\subset\C^n$ be an open set ($n\ge 1$) and let
$\Sigma\subset U\times\C^q$ for $q\ge2$ be a closed analytic
subset such that each fiber
$\Sigma_x =\{w\in \C^q\colon (x,w)\in \Sigma\}$ has complex
codimension at least two in $\C^{q}$ (it may be empty).
Assume that there exists a nonempty open set
$\Omega \subset \CP^{q-1}$ such that for
each $[v]\in \Omega$ the linear projection
$\wt \pi_v\colon U\times \C^q\to U\times \C^{q-1}$,
$$
    \wt \pi_v (x,w)=
      \bigl( x,\pi_v(w) \bigr)
      \quad (x\in U,\ w\in \C^q),               
$$
is proper when restricted to $\Sigma$. Then the submersion
$h\colon (U\times\C^q)\bs \Sigma \to U$, $h(x,w)=x$,
admits a spray.

Observe that condition (1.1) in part (a) of Theorem 1.6 implies that
for each vector $v\in \C^q$ sufficiently close to $e_q=(0,\ldots,0,1)$,
the projection $\pi_v$ is proper on $\Sigma$. Granted Lemma 7.1
we get a spray on $\C^q\bs \Sigma$, and hence part (a) of Theorem 1.6
follows from Theorem 1.2.

\demo Proof of Lemma 7.1:
We denote the coordinates on $U\times \C^q$ by $z=(x,w)$. Fix a point
$z_0=(x_0,w_0) \in (U\times\C^q)\bs \Sigma$. For each
line $[v]\in \Omega$ such that the affine complex
line $\{w_0+tv\colon t\in \C\} \subset \C^q$ does not intersect
$\Sigma_{x_0}$, the projection $\wt \pi_v$  satisfies
$\wt \pi_v(z_0)\notin \wt \pi_v(\Sigma)$.
By the codimension condition on $\Sigma$ this holds
for all directions $[v]$ outside a proper subvariety of $\Omega$.
For each such $[v]$ there exists a holomorphic function
$g\colon U\times \C^{q-1} \to \C$ which vanishes on the subvariety
$\wt \pi_v(\Sigma)$ and equals one at the point $\wt \pi_v(z_0)$.
The \hvf\ $V(z)= g\bigl( \wt\pi_v(z) \bigr) v$
(a shear field) is $\C$-complete and vertical on $U\times \C^q$,
with the flow $\theta_t(z)= z+ t g\bigl( \wt \pi_v(z)\bigr)v$\ \
($t\in \C$). By the choice of $g$ the field $V$ vanishes on
$\Sigma$, and hence its restriction to $(U\times\C^q)\bs\Sigma$
is also a complete field. We have $V(z_0)=v$.

This shows that the complete vertical fields on
$(U\times\C^q)\bs\Sigma$ generate the vertical tangent space
at each point. It remains to see that we can do the same with
finitely many fields. We may assume that $U$ is connected.
We begin by choosing complete vertical fields $V_1,\ldots, V_q$
as above which generate $VT_z(U\times\C^q)$ at one point of
$(U\times \C^q)\bs \Sigma$. Hence the same fields generate the tangent
space at each point outside a proper analytic subset
$A\subset (U\times \C^q)\bs \Sigma$. Let $A=\cup_j A_j$
be the (finite or countable) decomposition of $A$ into
irreducible components. Choose a point $z_j\in A_j$ in each
component and consider the set $\Omega_j\subset \Omega$ of
all complex directions $[v]\in \Omega$ for which $v$ belongs to
the linear span of the vectors $V_k(z_j)$ ($1\le k\le q$) or
$\wt \pi_v(z_j)\in \wt \pi_v(\Sigma)$. Clearly $\Omega_j$
is a proper analytic subset of $\Omega$ and hence $\cup_j \Omega_j$
is a set of the first category in $\Omega$. Choose
any direction $[v]\in  \Omega\bs \cup_j \Omega_j$. By construction we
then have $\wt \pi_v(z_j)\notin \wt \pi_v(\Sigma)$ for each $j$.
Let $g_1,\ldots,g_k$ be holomorphic functions in $U\times \C^{q-1}$
whose common zero set is precisely $\wt \pi_v(\Sigma)$.
If we add the corresponding complete vertical vector fields
$W_i(z)= g_i\bigl( \wt\pi_v(z) \bigr)v$ for $i=1,\ldots,k$
to the previous collection $V_1,\ldots,V_q$, we increase the
dimension of the linear span by at least one at each point $z_j$,
and hence at each point outside a proper subvariety of each
irreducible component $A_i$ of $A$. An induction on dimensions
of the span and of the exceptional set completes the proof
of Lemma 7.1.

\demo Proof of Theorem 1.6 (b):  According to Rosay and Rudin [RRu]
there is for each $q\in \N$ a discrete set $\Sigma\subset\C^q$
such that any entire holomorphic map $f\colon \C^n \to \C^q\bs \Sigma$
has complex rank at most $q-1$ at each point. The same is then
true for holomorphic maps $f\colon X\to \C^q\bs \Sigma$ from any Stein
manifold $X$ whose universal cover is a Euclidean space, such as
$X=(\C^*)^n$. We claim that any such map is holomorphically homotopic
to a constant map of $X$ into $\C^q\bs \Sigma$.
To see this, note that the rank condition on $f$ implies that
$f(X)$ is contained in a countable union
$A=\cup_{j=1}^\infty A_j \subset \C^q\bs \Sigma$, where each
$A_j$ is a compact subset contained in a local analytic set
of complex dimension $\le q-1$ in $\C^q$ [Chi].
We denote by $C_z(A_j)$ the real cone on $A_j$ with vertex at
$z\in \C^q$. Since $\Sigma$ is discrete and $C_z(A_j)$ has dimension
at most $2q-1$, we have $C_z(A_j) \cap\Sigma =\emptyset$
for an open and dense set of points $z\in \C^q$.
By Baire's theorem we can choose $z \in\C^q \bs\Sigma$ such
that the above holds for all $j$ and hence
$C_z(A) \cap \Sigma =\emptyset$. The homotopy
$f_t(x)=tz+(1-t)f(x)$ $(0\le t\le 1$) contracts the initial
map $f=f_0$ to the constant map $f_1(x)=z$ in $\C^q\bs \Sigma$,
establishing our claim.

On the other hand, when $n=2q-1$, there exist continuous
(even real-analytic) maps $f\colon X=(\C^*)^n \to \C^q\bs\Sigma$
which are not homotopic to constant: we contract
$(\C^*)^n$ onto the torus $T^n$ and embed $T^n$ as
a real hypersurface in $\C^q\bs \Sigma$ so that at least
one point of $\Sigma$ is contained in the bounded component
of $\C^q\bs T^n$. Thus the Oka principle fails for maps
$(\C^*)^{2q-1}\to \C^q\bs\Sigma$.

\demo Proof of Theorem 1.6 (c):
In [BFo] and [Fo1] it was shown that for any pair of integers
$1\le k<q$ there exist proper holomorphic embeddings $g\colon \C^k\to\C^q$
such that every entire map $f\colon \C^n \to D=\C^q\bs g(\C^k)$
has complex rank $<q-k$ at each point. The same proof as in case
(b) shows that, if $X$ is any Stein manifold which is universally
covered by a Euclidean space, any \holo\ map $X\to D=\C^q\bs g(\C^k)$
is contractible to a point in $D$, while for some such $X$ there
exist nontrivial \ra\ maps into $D$. For instance, choose a point
$z\in g(\C^k)$ and let $\Lambda\subset \C^q$ be the normal plane
to $g(\C^k)$ at $z$ (of complex dimension $q-k$). We can embed
the torus $T^n$, with $n=2(q-k)-1$, into $\Lambda$ so that $z$ is
contained in the bounded component of $\Lambda\bs T^n$ in $\Lambda$.
Since $T^n$ is a retraction of $(\C^*)^n$, we get a map
$(\C^*)^{2(q-k)-1} \to D$ which is not zero homotopic in $D$.
\endpr

%
%
\beginsection 8. Sections of vector bundles avoiding analytic subsets.

\demo Proof of Theorem 1.7:
If $X_0=\emptyset$, Theorem 1.7 follows from
Theorem 1.5 and Lemma 7.1 which may be seen as follows.
To apply Theorem 1.5 we need to show that for each $x\in X\bs K$
there is a \nbd\ $U\subset X$ of $x$ such that the submersion
$\pi^{-1}(U)\bs \Sigma \to U$ admits a spray.
Let $U$ and $\Phi$ be chosen as in Theorem 1.5.
Set $U'=\pi^{-1}(U)$ and
$\Sigma'=\Phi(\Sigma\cap U') \subset U\times\C^q$.
Condition (1.2) in Theorem 1.7 implies that for each
$v\in \C^q$ sufficiently close to $e_q=(0,\ldots,0,1)$,
the projection $\wt\pi_v\colon U\times \C^q\to U\times \C^{q-1}$,
given by $\wt \pi_v(x,w)=(x,\pi_v(w))$, is proper
on $\Sigma'$. By Lemma 7.1 the submersion
$(U\times \C^q)\bs\Sigma'\to U$ admits a spray $s'$ defined on a
trivial bundle over the given set. Then $s=\Phi^{-1} \circ s'$
is a spray associated to $U'\bs \Sigma \to U$ and hence Theorem 1.5
applies.

Interpolation along $X_0$ needs additional work. A direct
approach would require modifications in the approximation
theorems and in the gluing lemma in [FP1]. While it would be possible
and not overly difficult to carry out these modifications, we present
here an alternative approach which requires patching only on
subsets that do not intersect $X_0\cup K$.

The section $f_0 \colon X\to V$ is assumed to be \holo\ in
a \nbd\ of $X_0\cup K$. This set has a basis of Stein \nbd s.

%
%
\proclaim 8.1 Lemma: {\rm (Notation as above.)}
We can write $f_0$ in the form
$$ f_0= \phi+ \sum_{j=1}^m h_j\, g^0_j,         \eqno(8.1) $$
where $\phi\colon X\to V$ is global \holo\ section of $V$,
$g^0_j\colon X\to V$ are continuous sections which are \holo\
in an open set $U_0\supset X_0\cup K$, and the functions
$h_1,\ldots,h_m \in {\cH}(X)$ vanish to order $k$ on $X_0$
and satisfy $X_0=\{x\in X\colon h_j(x)=0,\ 1\le j\le m\}$.
(The graph of $\phi$ may intersect $\Sigma$ outside a \nbd\ of
$X_0$.)

\demo Proof: This is a straightforward application of the
Oka-Cartan theory, but we include the proof for completeness.
We begin by choosing finitely many functions
$h_1,\ldots,h_m \in {\cH}(X)$ with the stated properties.
We denote by ${\cO}= {\cO}_X$ the sheaf of germs of \holo\ functions
on $X$ and by ${\cV}$ the sheaf of germs of \holo\ sections
of $V\to X$. Furthermore, let ${\cJ} \subset {\cO}$ denote
the sheaf of ideals in ${\cO}$ generated by the functions
$h_1,\ldots,h_m$, and let ${\cV}_0={\cJ}{\cV} \subset {\cV}$ be
the sheaf of germs of sections of $V$ which can be locally written
as $\sum_{j=1}^m h_j g_j$ for some $g_j \in {\cV}$,
$1\le j\le m$. We have a short exact sequence of
coherent analytic sheaves on $X$:
$$ 0\to {\cV}_0 \hra {\cV} \to {\cK} \to 0. $$
The quotient sheaf $\cK={\cV}/{\cV}_0$ is trivial on $X\bs X_0$
by the choice of the $h_j$'s. Hence the section $f_0 \colon X\to V$,
which is \holo\ in an open set containing $X_0$, determines a global
\holo\ section $\wt f_0\colon X\to {\cK}$. By Cartans's Theorem B
we have $H^1(X;{\cV}_0)=0$, and hence there is a \holo\ section
$\phi\colon X\to {\cV}$ whose image in $\cK$ equals
$\wt f_0$. We identify $\phi$ with the corresponding section
of $V\to X$; then $f_0-\phi$ is a \holo\ section of $V\to X$
over an open set $U_1\supset X_0\cup K$ (which we may take to be
Stein), and by construction $f_0-\phi$ vanishes to order $k$ on
$X_0$.

Consider the sheaf epimorphism ${\cV}^m\to {\cV}_0$,
$(g_1,\ldots,g_m)\to \sum_{j=1}^m h_j g_j$.
By Cartan's Theorem B [GRo] we can lift each global \holo\
section of ${\cV}_0$ over the Stein manifold $U_1$ to a \holo\
section of ${\cV}^m \to U_1$. Applying this to $f_0-\phi$ we
obtain \holo\ sections $g^0_j\colon U_1\to V$ $(1\le j\le m)$
satisfying (8.1) over $U_1$.

We denote by $\wh V=V^m\to X$ the direct (Whitney) sum of $m$
copies of the bundle $V\to X$. Then
$G_0=(g^0_1,\ldots,g^0_m) \colon U_1\to \wh V$ is a \holo\
section of $\wh V \to U_1$. Using a partition of unity we can
modify $G_0$ outside a smaller \nbd\ $U_0\subset U_1$ of $X_0\cup K$
(without changing it on $U_0$) and extend it continuously to $X$
such that (8.1) holds everywhere.
\endpr

Denote by $\Theta\colon \wh V \to V$ the map
$$
   \Theta(x;v_1,\ldots,v_m) =
   \bigl( x; \phi(x)+ \sum_{j=1}^m h_j(x) v_j \bigr), \eqno(8.2)
$$
where  $x\in X$ and $v_1,\ldots,v_m \in V_x$. Clearly $\Theta$ is a
submersion over $X\bs X_0$ and is degenerate over $X_0$.
Set $\wh \Sigma=\Theta^{-1}(\bar \Sigma) \subset \wh V$.
Then the graph of a section $G=(g_1,\ldots,g_m) \colon X \to \wh V$
avoids $\wh \Sigma$ \iff\ the graph of the associated section
$f=\Theta\circ G=\phi+\sum_{j=1}^m h_jg_j \colon X\to V$
avoids $\bar \Sigma$. Notice that every section of this type
agrees with $f_0$ to order $k$ along $X_0$. By construction
we have $G_0(X) \cap \wh \Sigma =\emptyset$, and $G_0$ is \holo\
in $U_0$. To complete the proof of Theorem 1.7 we need the following.

%
%
\proclaim 8.2 Proposition: {\rm (Notation as above.)}
Let $G_0 \colon X\to \wh V\bs \wh\Sigma$ be a continuous
section which is \holo\ in an open set $U_0\supset X_0\cup K$,
where $K \subset X$ is \hc{X}.
For each compact \hc{X}\ subset $L\subset X$ containing $K$
there are an open set $U'_0$, with $X_0\cup K \subset U'_0\subset U_0$,
and a homotopy $G_t \colon X\to \wh V \bs\wh \Sigma$ $(t\in [0,1])$
of continuous sections which are \holo\ in $U'_0$ such that
$G_t$ approximates $G_0$ uniformly on $K$ for each $t\in [0,1]$
and the section $G_1$ is \holo\ in an open set $W_0\supset X_0\cup L$.

Granted Proposition 8.2 we can complete the proof of
theorem 1.7 as follows. We exhaust $X$ by a sequence
$\{L_j\colon j\in \Z_+\}$ of compact \hc{X}\ sets
such that $L_0=K$ and $X=\cup_{j=0}^\infty L_j$.
Applying Proposition 8.2 inductively on each pair
$(L_j,L_{j+1})$ we obtain a sequence
of sections $G_j\colon X\to \wh V\bs \wh\Sigma$ $(j\in \Z_+)$
such that for each $j$, $G_{j+1}$ is \holo\ in a \nbd\
of $X_0\cup L_{j+1}$ and is homotopic to $G_j$ by a
homotopy $G_t\colon X\to \wh V\bs \wh \Sigma$
$(t\in [j,j+1])$ which is \holo\ in a \nbd\ $X_0\cup L_j$
and which approximates $G_j$ uniformly on $L_j$.
If these approximations are sufficiently close, the sequence
$G_j$ converges to a \holo\ section
$G_\infty =\lim_{j\to\infty} G_j \colon X\to \wh V\bs \wh \Sigma$.
By reparametrizing the homotopy $G_t$ $(t\in [0,+\infty))$
as in the proof of Theorem 1.5 (sect.\ 6) we get a homotopy
$G_t$ $(t\in [0,1])$ from $G_0$ to $G_1=G_\infty$.
The associated homotopy
$f_t=\Theta\circ G_t \colon X \to V$ $(t\in[0,1])$ then
satisfies Theorem 1.7.
\endpr

To prove Proposition 8.2 we shall carry out the modification
procedure, described in the proof of Theorem 1.5, so that
we only glue sections over Cartan pairs $(A,B)$ in $X$ for which
$B\cap (X_0\cup K)=\emptyset$. We need the following lemmas.

%
%
\proclaim 8.3 Lemma:
The submersion $\wh\pi \colon \wh V\bs\wh \Sigma \to X$ admits a
spray over a \nbd\ of any point $x\in X'=X\bs (X_0\cup K)$.

\demo Proof: Since the map $\Theta$ (8.2) is an affine linear
submersion of \hvb s over $X\bs X_0$, it is locally (over small
sets $U\subset X\bs X_0$) equivalent to the projection $(x;v,w)\to (x;v)$
of trivial bundles. In such coordinates on $\wh V|_U$ resp.\ $V|_U$
the set $\wh \Sigma \cap \wh\pi^{-1}(U)$ is defined by the same
equations as $\Sigma \cap \pi^{-1}(U)$, and the additional $(m-1)q$
fiber coordinates are not present in these equations. It is now
immediate that the validity of (1.2) for $\Sigma\cap U'$ implies the
analogous condition for $\wh \Sigma\cap \wh U$. Lemma 7.1 shows
as before that the submersion $\wh \pi \colon \wh V\bs \wh \Sigma \to X$
admits a spray over a \nbd\ of any point $x\in X'$.
\endpr

%
%
%
\medskip\ni \bf 8.4 Lemma.  \sl
Let $X$ be a Stein manifold, $X_0$ a closed analytic subvariety of
$X$ and $K, L \subset X$ a pair of compact \hc{X}\ subsets,
with $K\subset {\rm Int}L$. Let $\cU=\{U_j\}_{j=0}^\infty$ be an
open covering of $X$ such that $X_0\cup K \subset U_0$ and
$(X_0\cup K)\cap U_j =\emptyset$ for each $j\ge 1$.
Then there is a Cartan string $(A_0,A_1,\ldots,A_n)$ in $X$
such that $A_0$ is \hc{X}\ and the following hold:

\item{(i)}   $K\cup (X_0\cap L) \subset A_0 \subset U_0$;

\item{(ii)}  for $j=1,2,\ldots,n$ we have
$A_j \cap (X_0 \cup K)=\emptyset$ and $A_j \subset U_k$
for some $k=k(j)\ge 1$;

\item{(iii)} $L=\bigcup_{0\le j\le n} A_j$.
\medskip\rm

\demo Proof:  Choose a compact \hc{X}\ set $K'\subset X$ with
$L\subset {\rm Int}K'$. Then $S=(X_0\cap K')\cup K$ is also \hc{X}.
By Theorem 4.6 there exists a Cartan string $(A'_0,A'_1,\ldots,A'_r)$
in $X$ satisfying the following:
\item{(a)}  $S\subset A'_0 \subset U_0$;
\item{(b)}  if $1\le j\le r$ then $A'_j \cap S=\emptyset$;
\item{(c)}  if $A_j\cap X_0=\emptyset$ then $A'_j\subset U_{k(j)}$
for some $k(j)\ge 1$;
\item{(d)}  if $A'_j \cap X_0 \ne\emptyset$ for some
$1\le j\le r$ then $A'_j\cap L=\emptyset$;
\item{(e)}  $L\subset \bigcup_{0\le j\le r} A'_j$.

Let $A_0,\ldots,A_n$ denote the nonempty sets in the string
$A'_j\cap L$, $0\le j\le r$ (in the given order). Since $L$ is \hc{X},
Proposition 4.3 implies that $(A_0,A_1,\ldots,A_n)$ is a Cartan string
in $X$. It is clear that (i) and (iii) in Lemma 8.3 hold, and (ii) holds
because no set $A'_j$ for $j\ge 1$ intersects both $X_0$ and $L$
at the same time according to (d). This proves Lemma 8.3.
\endpr

%
%
\demo Proof of Proposition 8.2:
By Lemma 8.3 there is an open covering ${\cal U }= \{U_j\}_{j=0}^\infty$
of $X$ such that the submersion $\wh\pi\colon \wh V\bs \wh \Sigma\to X$
admits a spray over each set $U_j$ for $j\ge 1$. Furthermore, we may
choose the sets in the covering such that Proposition 4.7 applies to
$G_0$, i.e., we can modify $G_0$ into a \holo\ $\cK(\cU)$-complex.

Let $\cA=(A_0,A_1,\ldots,A_n)$ be a Cartan string provided by
Lemma 8.4, subordinate to ${\cal U}$ and satisfying
$L=\cup_{j=0}^n A_j$. Applying Proposition 4.7 we deform $G_0$
over a \nbd\ of $L$ into a holomorphic $\cK(\cA)$-complex
$H_*$. Since $G_0$ is holomorphic on $U_0 \supset A_0$,
we may (and do) take the section $H_{(0)}$ in the
complex $H_*$ to be $G_0$, restricted to a \nbd\ of $A_0$.

By the process described in section 5 (see especially the
final remark at the end of sect.\ 5) we can modify the complex
$H_*$ in a finite sequence of steps into a holomorphic section
$H_1\colon W\to \wh V \bs\wh \Sigma$ over an open
set $W\supset L$ such that $H_1$ approximates $G_0$
uniformly on a \nbd\ $U_1\supset A_0$. In addition,
shrinking $U_1$ around $A_0$ if necessary, the construction
in sect.\ 5 gives a homotopy of sections
$H_t \colon W \to\wh V\bs \wh \Sigma$ $(t\in [0,1])$
which are \holo\ in $U_1$ and approximate $G_0$ there,
with $H_0=G_0$.

Since $L$ is \hc{X}, we may approximate $H_1$ (which is \holo\ in
$W$) uniformly on $L$ by a global \holo\ section $G\colon X\to V$.
For $t\in [1,2]$ set $H_t=(2-t)H_1+(t-1)G \colon X\to \wh V$. At
$t=1$ this coincides with the section $H_1$ defined earlier. We
also have $H_2=G$, and for each $t\in [1,2]$ the section $H_t$ is
\holo\ in $W$ and it approximates $G_0$ in $U_1\supset A_0$. For
convenience of notation we replace $t \in [0,2]$ by $t/2\in [0,1]$,
thus reparametrizing $\{H_t\colon t\in [0,2]\}$ to the interval
$t\in [0,1]$.

By shrinking the set $U_0\supset X_0\cup A_0$ we may assume that
$U_0$ is Stein. We can now approximate the homotopy
$H_t \colon U_1\to \wh V$, uniformly on $A_0 \subset U_1$,
by a \holo\ homotopy
$\wt H_t\colon U_0 \to \wh V$ $(t\in [0,1])$
satisfying $\wt H_0=G_0$ and $\wt H_1=G$. This
can be done for instance by applying the h-Runge approximation
theorem to the family $H_t$ twice.
First we apply it with the pair $A_0\subset U_0$ and the initial
section $G_0$ to obtain a family $H^{(1)}_t\colon U_0\to \wh V$
$(t\in [0,1])$ with $H^{(1)}_0=G_0$; the second time we apply it
with the pair $A_0\subset X$ and the `initial' section $H_1$ to get
a family $H^{(2)}_t\colon X\to \wh V$ $(t\in [0,1])$ with
$H^{(2)}_1=G$. All sections $H^{(1)}_t$ and $H^{(2)}_t$ approximate
$H_t$, and hence $G_0$, uniformly on $A_0$. Finally we take
$$ \wt H_t=(1-t)H^{(1)}_t + t H^{(2)}_t \qquad (t\in [0,1])$$
(restricted to $U_0$ when $t<1$).

The homotopy $\wt H_t$ is \holo\ in $U_0$.
Recall that, over $A_0$, both $H_t$ and $\wt H_t$
approximate $G_0$ as close as desired, and $\wh \Sigma$ has
no points over $X_0$. It follows that, if the approximations
are sufficiently close, there is an open set $U'_1\subset X$,
with $X_0\cup A_0 \subset U'_1\subset U_0$, such that for any
$x\in U'_1$ and $t,\tau \in [0,1]$ we have
$$ \tau \wt H_t(x) + (1-\tau) H_t(x)\notin\wh\Sigma.  \eqno(8.3)$$
Furthermore the graph of $H_1=G$ avoids $\wh \Sigma$ over an open
set $W'_0\supset X_0\cup L$. Choose a smaller open set
$U'_0\subset X$, with $X_0\cup A_0\subset U'_0$ and
$\bar U'_0\subset U'_1$, and choose a smooth function
$\tau \colon X\to [0,1]$ satisfying $\tau =1$ on
$\bar U'_0$ and $\supp \tau \subset U'_1$. Define
a new homotopy by
$$
    G'_t(x)= \tau(x) \wt H_t(x) + (1-\tau(x)) H_t(x) \qquad
    (x\in U'_0\cup W,\ t\in [0,1]).
$$
The first term $\wt H_t$ is only defined on $U'_1$, but
since $\tau=0$ on $X\bs U'_1$, we can extend this terms to $X$.
Likewise $H_t$ is only defined on $W$, but $1-\tau(x)=0$ on
$U'_0\supset X_0$ and hence  the second term extends to $U'_0 \cup W$.
Thus $G'_t$ is defined on $U'_0\cup W$ and it satisfies
the following:

\item{(i)} $G'_0=G_0$ (since $\wt H_0=H_0=G_0$),
and likewise $G'_1=G$ (since $\wt H_1=H_1=G$);

\item{(ii)} for each $t\in [0,1]$ the section $G'_t$
is \holo\ in $U'_0$;

\item{(iii)} shrinking $W\supset L$ if necessary we have
$G'_t(x) \notin \wh\Sigma$ for any $x\in U'_0\cup W$
and $t\in [0,1]$. For $x\in U'_0 \cup (W\cap U'_1)$ this holds
by (8.3), while for $x\in W\bs U'_1$ we have $\tau(x)=0$ and
hence $G'_t(x)= H_t(x) \notin \wh\Sigma$.

It remains to extend the homotopy $G'_t$ to $X$.
Choose an open set $W_0\subset X$, with $X_0\cup L\subset W_0$
and $\bar W_0\subset U'_0\cup W$, and choose a smooth function
$\chi \colon X\to [0,1]$ which equals one on $W_0$ and zero on
$X\bs (U'_0\cup W)$. The homotopy
$$
    G_t(x)= G'_{\chi(x)t}(x) \qquad (x\in X,\ t\in [0,1])
$$
then satisfies Proposition 8.2.
\endpr

\demo Proof of Theorem 1.9: The scheme is the same as in
Theorem 1.7. Using the notation established above, we construct a
section $G_0$ of the submersion
$\wh \pi\colon \wh V\bs\wh \Sigma \to X$ such that
$\Theta\circ G_0=f_0$ is the given section
of $V\bs \Sigma\to X$, where $\Theta$ is the map (8.2).
Our goal is to establish Proposition 8.2 in this situation, and
for this we must show that Lemma 8.3 holds in the current setting,
i.e., our submersion admits a spray over a small \nbd\ $U$ of any
point $x\in X\bs X_0$.

Since $\wh\pi(\wh\Sigma)=X_1\bs X_0$, we only need to consider
points $x\in X_1\bs X_0$. For $U\subset X$ we set
$U'=\pi^{-1}(U)\subset V$ and
$\wh U=\wh \pi^{-1}(U) \subset \wh V$.
In the proof of Lemma 8.2 we saw that, over a small \nbd\
$U$ of $x$, there is a fiber preserving biholomorphic map
$\Phi\colon \wh U \to U'\times \C^{(m-1)q}$ which maps
$\wh \Sigma \cap \wh U$ onto $(\Sigma \cap U') \times \C^{(m-1)q}$.
By the assumption, if $U$ is chosen sufficiently small,
there is a \holo\ action of a complex Lie group $G$ on $U'$
which preserves the fibers of $\pi$ and acts transitively on
$V_y\bs \Sigma_y$ for any $y\in U\cap X_1$. The Lie group
$G_1=G\times \C^{(m-1)q}$ then acts on $U'\times \C^{(m-1)q}$ by
$$ (g,w)\cdotp (z,w')=(g\cdotp z,w+w') \qquad
    (g\in G,\ z\in U',\ w,w'\in \C^{(m-1)q}).
$$
Conjugating by $\Phi$ we get an associated action of $G_1$
on $\wh U$ which preserves the fibers of $\wh \pi$ and
acts transitively on $\wh V_y\bs\wh\Sigma_y$ for
$y\in U\cap X_1$.

Let ${\bf g_1}$ denote the Lie algebra of $G_1$. The holomorphic
map $s_1\colon E_1=\wh U\times {\bf g_1} \to \wh U$, defined by
$s_1(z;g)=\exp(g)\cdotp z$ for $z\in \wh U$ and $g\in {\bf g_1}$,
satisfies all conditions for a spray over $\wh U\bs \wh \Sigma$
(def.\ 1.1), except that its vertical derivative at the zero section
need not be surjective at points $z\in \wh U$ for which
$\wh\pi(z) \in U\bs X_1$. This is easily corrected as
follows. Choose \holo\ functions $b_1,\ldots, b_k$ in $U$ such that
$X_1\cap U=\{x\in U\colon b_j(x)=0,\ 1\le j\le k\}$.
Also choose $\C$-complete vertical \hvf s $V_1,\ldots,V_{mq}$
on $\wh U$ which generate the vertical tangent space at each
point. (Since $\wh U\cong U\times \C^{mq}$, we may take the
constant fields in the $mq$ coordinate directions on $\C^{mq}$.)
Let $\{\phi_l^t\colon 1\le l\le N\}$ be the flows of the
$N=kmq$ complete fields $(b_j\circ \wh\pi)V_i$ for $1\le i\le mq$ and
$1\le j\le k$. Set $E=\wh U\times ({\bf g_1} \oplus \C^N) \to \wh U$
and let $s\colon E\to \wh U$ be the map
$$
    s(z;g,t_1,\ldots,t_N) =
     \phi_1^{t_1}\circ\cdots\circ\phi_N^{t_N}(\exp(g)\cdotp z), \eqno(8.2)
$$
where $z\in\wh U$, $g\in {\bf g_1}$ and $t_j\in \C$ for $1\le j\le N$.
It is immediate that $s$, restricted to the part of the bundle
$E$ over $\wh U\bs \wh \Sigma$, is a spray on
$\wh\pi\colon \wh U\bs \wh \Sigma \to U$.
\endpr

%
%
%
\medskip \ni \it Acknowledgements. \rm
The first author acknowledges partial support by the
NSF, by the Vilas Foundation at the University of Wisconsin--Madison,
and by the Ministry of Science and Technology of the Republic of Slovenia.
The second author was supported by the Ministry of Education
of the Republic of Slovenia. We thank the participants in
the complex analysis seminar at the IMFM in Ljubljana for the
opportunity to explain the details of this proof in a series of
lectures in Fall and Winter of 1998. The first author did a
part of the work while visiting Universit\'e Paul Sabatier in 
Toulouse on several occasions in 1998-99. He wishes to thank the 
colleagues at this institution for invitations, interesting discussions 
and their kind hospitality. We also thank J.\ Globevnik,
G.\ M.\ Henkin, J.\ Leiterer and R.\ Narasimhan for their interest 
in these developments.

\medskip\ni \it Addendum. \rm
This paper was completed in February 1999. In the mean time its 
sequel [FP2], which depends on the methods developed here, 
has been published. [FP2] contains the Oka principle with interpolation 
of sections on closed complex subvarieties of the base manifold. 
Recently F.\ L\'arusson [Lar] gave a homotopy theoretic proof of 
the Oka principle, granted the analytic results on approximation 
and patching of \holo\ sections proved in [FP1]. The papers 
[Pre1, Pre2] contain new applications of the Oka principle to 
the embedding problem.  In [Fo2] and [Fo3] these methods were used
to study removability of intersections of \holo\ maps from Stein 
manifolds with complex subvarieties of the target manifold. 
The paper [Fo4] contains a version of the Oka principle for 
multi-valued sections of certain ramified holomorphic maps over 
a Stein base. In [Fo5] the Oka principle is extended to a 
wider class of target manifolds and submersions.

%
%
%
%
\medskip\ni\bf References. 
\rm\medskip

\ii{[BFo]} G.\ Buzzard, J.\ E.\ Forn\ae ss:
An embedding of $\C$ in $\C^2$ with hyperbolic complement.
Math.\ Ann.\ {\bf 306}, 539--546 (1996).

\ii{[Car]} H.\ Cartan: Espaces fibr\'es analytiques.
Symposium Internat.\ de topologia algebraica, Mexico, 97--121 (1958).
(Also in Oeuvres, vol.\ 2, Springer, New York, 1979.)

\ii{[Chi]} E.\ M.\ Chirka: Complex analytic sets.
Mathematics and its Applications (Soviet Series), {\bf 46},
Kluwer, Dordrecht, 1989.

\ii{[EGr]} Y.\ Eliashberg, M.\ Gromov: Embeddings of Stein manifolds.
Ann.\ Math.\ {\bf 136}, 123--135 (1992).

\ii{[Fo1]} F.\ Forstneri\v c:
Interpolation by holomorphic automorphisms and embeddings in $\C^n$.
J.\ Geom.\ Anal.\ {\bf 9}, no.1, (1999) 93-118.

\ii{[Fo2]} F.\ Forstneri\v c: On complete intersections.
Ann.\ Inst.\ Fourier {\bf 51} (2001), 497--512.

\ii{[Fo3]} F.\ Forstneri\v c:
The Oka principle, lifting of holomorphic maps and
removability of intersections. 
Proc.\ of Hayama Symposium on Several Complex Variables 2000,
pp.\ 49--59, Japan, 2001.

\ii{[Fo4]} F.\ Forstneri\v c:
The Oka principle for multivalued sections of ramified mappings.
Forum Math., to appear.

\ii{[Fo5]} F.\ Forstneri\v c:
The Oka principle for sections of subelliptic submersions. 
Math.\ Z., to appear.

\ii{[FP1]} F.\ Forstneri\v c and J.\ Prezelj:
Oka's principle for sections of holomorphic fiber bundles with sprays.
Math.\ Ann.\ {\bf 317} (2000), 117--154.

\ii{[FP2]} F.\ Forstneri\v c and J.\ Prezelj:
Extending holomorphic sections from complex subvarieties.
Math.\ Z.\ {\bf 236} (2001), 43--68.

\ii{[Gra]} H.\ Grauert:
Holomorphe Funktionen mit Werten in komplexen Lieschen Gruppen.
Math.\ Ann.\ {\bf 133}, 450--472 (1957).

\ii{[Gr1]} M.\ Gromov:
Oka's principle for holomorphic sections of elliptic bundles.
J.\ Amer.\ Math.\ Soc.\ {\bf 2}, 851-897 (1989).

\ii{[Gr2]} M.\ Gromov: Partial differential relations.
Ergebnisse der Mathematik und ihrer Grenzgebiete (3) {\bf 9},
Springer, Berlin--New York, 1986.

\ii{[GRo]} C.\ Gunning, H.\ Rossi:
Analytic functions of several complex variables.
Prentice--Hall, Englewood Cliffs, 1965.

\ii{[HL]} G.\ Henkin, J.\ Leiterer:
The Oka-Grauert principle without induction over the basis dimension.
Math.\ Ann.\ {\bf 311}, 71--93 (1998).

\ii{[H\"or]} L.\ H\"ormander:
An Introduction to Complex Analysis in Several Variables, 3rd ed.
North Holland, Amsterdam, 1990.

\ii{[HW]} W.\ Hurewicz, H.\ Wallman: Dimension Theory.
Princeton Mathematical Series, {\bf 4}, Princeton University Press,
Princeton, 1941.

\ii{[Lar]} F.\ L\'arusson:
Excision for simplicial sheaves on the Stein site
and Gromov's Oka principle. Preprint, December 2000.

\ii{[Oka]} K.\ Oka: Sur les fonctions des plusieurs variables. III:
Deuxi\`eme probl\`eme de Cousin.
J.\ Sc.\ Hiroshima Univ.\ {\bf 9} (1939), 7--19.

\ii{[Pre1]} J.\ Prezelj: Interpolation of embeddings of Stein
manifolds on discrete sets. Pre\-print, 1999.

\ii{[Pre2]} J.\ Prezelj: Weakly regular embeddings of Stein spaces.
Preprint, 2001.

\ii{[RRu]} J.-P.\ Rosay, W.\ Rudin:
Holomorphic maps from $\C^n$ to $\C^n$.
Trans.\ Amer.\ Math.\ Soc.\ {\bf 310}, 47--86 (1988).

\ii{[Sch]} J.\ Sch\"urmann:
Embeddings of Stein spaces into affine spaces of minimal dimension.
Math.\ Ann.\ {\bf 307}, 381--399 (1997).

\ii{[Spr]} D.\ Spring: Convex integration theory 
(Solutions to the $h$-principle in geometry and topology).
Monographs in Mathematics, {\bf 92}, Birkh\"auser, Basel, 1998.

\ii{[Win]} J.\ Winkelmann:
The Oka-principle for mappings between Riemann surfaces.
Enseign.\ Math.\ {\bf 39} (1993), 143--151.

\bigskip
\ni F.\ Forstneri\v c, J.\ Prezelj
\par \ni
Institute of Mathematics, Physics and Mechanics,
University of Ljubljana, Jadranska 19, 
\ni 1000 Ljubljana, Slovenia (e-mail: Franc.Forstneric@fmf.uni-lj.si)

\bye